\documentclass[preprint]{elsarticle}
\usepackage[bookmarks=true,colorlinks=true,linkcolor=blue]{hyperref}

\usepackage{amsmath,amssymb,mathabx,mathtools}
\usepackage{graphicx,esint,ulem}

\biboptions{sort&compress}
\usepackage[usenames,dvipsnames,svgnames,table]{xcolor}

\usepackage{color}
\definecolor{jwbGreen}{rgb}{0, .6, 0}
\definecolor{jbaPurple}{HTML}{6600FF}
\definecolor{purple}{rgb}{.7, 0., .8}

\newcommand{\blue}{\color{blue}}

\usepackage{calc}
\usepackage[margin=1.in]{geometry}

\usepackage{fancyvrb}

\usepackage{empheq}
\usepackage[most]{tcolorbox}

\newtcbox{\mymath}[1][]{%
    nobeforeafter, math upper, tcbox raise base,
    enhanced, colframe=blue!60!black,
    colback=blue!20, boxrule=1pt,
    #1}

\usepackage{algorithm,algpseudocode}
\usepackage{algorithmicx}
\algrenewcommand\alglinenumber[1]{\footnotesize #1:} 
\newcommand{\algFontSize}{\footnotesize}

\usepackage{tcolorbox}
\usepackage[english]{babel}
\usepackage{amsmath}
\usepackage{amsfonts}
\usepackage{fullpage}
\usepackage{graphicx}
\usepackage{setspace}
\usepackage{float}

\usepackage{tabu}

\usepackage{multirow}

\usepackage{mdframed}

\allowdisplaybreaks

\usepackage{tikz,pgfplots}

\newtheorem{thm}{Theorem}[section]

\newcommand{\citeCount}[1]{}

\newcommand{\bogus}[1]{{}}

\newcommand{\mni}{\medskip\noindent}

\newcommand{\CR}{{\rm CR}}

\newlength{\ycbTop}
\newlength{\ycbMid}%

\newcommand{\p}{\partial}

\newcommand{\Np}{{N_p}}

\newcommand{\f}[2]{\frac{#1}{#2}}

\def\ba#1\ea{\begin{align}#1\end{align}}

\def\bas#1\eas{\begin{align*}#1\end{align*}}

\def\bat#1\eat{\begin{alignat}{3}#1\end{alignat}}

\def\bats#1\eats{\begin{alignat*}{3}#1\end{alignat*}}

\newcommand{\bse}{\begin{subequations}}
\newcommand{\ese}{\end{subequations}}

\newcommand{\Dzt}{D_{0t}}
\newcommand{\Dpt}{D_{+t}}
\newcommand{\Dmt}{D_{-t}}
\newcommand{\dt}{\Delta t}
\newcommand{\dx}{\Delta x}

\newcommand{\dr}{{\Delta r}}

\newcommand{\eqdef}{\overset{{\rm def}}{=}}

\newcommand{\bv}{\mathbf{ b}}

\newcommand{\gv}{\mathbf{ g}}

\newcommand{\jv}{\mathbf{ j}}

\newcommand{\pv}{\mathbf{ p}}

\newcommand{\rv}{\mathbf{ r}}
\newcommand{\sv}{\mathbf{ s}}
\newcommand{\tv}{\mathbf{ t}}

\newcommand{\vv}{\mathbf{ v}}
\newcommand{\wv}{\mathbf{ w}}
\newcommand{\xv}{\mathbf{ x}}
\newcommand{\yv}{\mathbf{ y}}
\newcommand{\zv}{\mathbf{ z}}

\newcommand{\Bv}{\mathbf{ B}}

\newcommand{\Gv}{\mathbf{ G}}

\newcommand{\Uv}{\mathbf{ U}}
\newcommand{\Vv}{\mathbf{ V}}

\newcommand{\Zv}{\mathbf{ Z}}

\newcommand{\half}{\f{1}{2}}

\newcommand{\Real}{{\mathbb R}}

\newcommand{\zerov}{\mathbf{0}}

\newcommand{\Ac}{{\mathcal A}}
\newcommand{\Bc}{{\mathcal B}}

\newcommand{\Dc}{{\mathcal D}}

\newcommand{\Gc}{{\mathcal G}}

\newcommand{\Kc}{{\mathcal K}}

\newcommand{\etav}{\boldsymbol{\eta}}
\newcommand{\xiv}{{\boldsymbol{\xi}}}
\newcommand{\muv}{{\boldsymbol{\mu}}}
\newcommand{\Phiv}{{\boldsymbol{\Phi}}}

\newcommand{\eps}{\epsilon}

\newcommand{\ds}{{\Delta s}}

\newcommand{\ssf}{\scriptscriptstyle}

\newcommand{\Tbar}{{\widebar T}}
\newcommand{\Tf}{\Tbar}

\newcommand{\sinc}{{\rm sinc}}

\newcommand{\omegaTilde}{\tilde{\omega}}

\newcommand{\Gcd}{\Gc_d}
\newcommand{\ACR}{{\rm ACR}}

\newcommand{\PPW}{{\rm PPW}}

\newcommand{\alphag}{a_{g}}
\newcommand{\betag}{b_{g}}

\newcommand{\maxit}{{\rm maxit}}

\DeclareMathOperator*{\argmin}{arg\,min}
\newcommand*{\vertbar}{\rule[-1ex]{0.5pt}{2.5ex}}

\newcommand{\kappaA}{\kappa_{\ssf A}}
\newcommand{\DeflationSet}{\Dc} 

\newcommand{\rvHat}{\hat{\rv}}
\newcommand{\rvTilde}{\tilde{\rv}}
\newcommand{\rhoOld}{\rho_{\rm old}}
\newcommand{\Lph}{L_{p,h}}

\newcommand{\Na}{N_{\DeflationSet}}
\newcommand{\Ndeflate}{\Na}

\newcommand{\kappaD}{\kappa_{\ssf D}}

\newcommand{\omegaN}{\omega_n}

\newcommand{\domainSize}{\mathcal{D}}

\newcommand{\Nwseig}{N_{\ssf EW}}

\newlength{\tfwidth}
\newlength{\tfheight}
\newlength{\tfxa}
\newlength{\tfxb}
\newlength{\tfya}
\newlength{\tfyb}
\newcommand{\trimFigWithBox}[6]{%
\setlength\fboxsep{0pt}%
\setlength\fboxrule{1.0pt}
\fbox{\includegraphics[width=#2, clip, trim=#3 #4 #5 #6]{#1}}%
}
\newcommand{\trimFigNoBox}[6]{%
\setlength\fboxsep{1pt}
\setlength\fboxrule{0.0pt}
\fbox{\includegraphics[width=#2, clip, trim=#3 #4 #5 #6]{#1}}%
}
\newcommand{\trimFigHeightWithBox}[6]{%
\setlength\fboxsep{0pt}%
\setlength\fboxrule{1.0pt}
\fbox{\includegraphics[height=#2, clip, trim=#3 #4 #5 #6]{#1}}%
}
\newcommand{\trimFigHeightNoBox}[6]{%
\setlength\fboxsep{1pt}
\setlength\fboxrule{0.0pt}
\fbox{\includegraphics[height=#2, clip, trim=#3 #4 #5 #6]{#1}}%
}

\newsavebox\figBox

\newcommand{\trimw}[6]{%
\sbox\figBox{\includegraphics{#1}}
\setlength{\tfwidth}{\the\wd\figBox}
\setlength{\tfheight}{\the\ht\figBox}
\setlength{\tfxa}{\tfwidth*\real{#3}}%
\setlength{\tfxb}{\tfwidth*\real{#4}}%
\setlength{\tfya}{\tfheight*\real{#5}}%
\setlength{\tfyb}{\tfheight*\real{#6}}%
\trimFigNoBox{#1}{#2}{\tfxa}{\tfya}{\tfxb}{\tfyb}%
}

\newcommand{\trimwb}[6]{%

\sbox\figBox{\includegraphics{#1}}
\setlength{\tfwidth}{\the\wd\figBox}
\setlength{\tfheight}{\the\ht\figBox}
\setlength{\tfxa}{\tfwidth*\real{#3}}%
\setlength{\tfxb}{\tfwidth*\real{#4}}%
\setlength{\tfya}{\tfheight*\real{#5}}%
\setlength{\tfyb}{\tfheight*\real{#6}}%
\trimFigWithBox{#1}{#2}{\tfxa}{\tfya}{\tfxb}{\tfyb}%
}

\newcommand{\trimh}[6]{%
\sbox\figBox{\includegraphics{#1}}
\setlength{\tfwidth}{\the\wd\figBox}
\setlength{\tfheight}{\the\ht\figBox}
\setlength{\tfxa}{\tfwidth*\real{#3}}%
\setlength{\tfxb}{\tfwidth*\real{#4}}%
\setlength{\tfya}{\tfheight*\real{#5}}%
\setlength{\tfyb}{\tfheight*\real{#6}}%
\trimFigHeightNoBox{#1}{#2}{\tfxa}{\tfya}{\tfxb}{\tfyb}%
}

\newcommand{\trimhb}[6]{%

\sbox\figBox{\includegraphics{#1}}
\setlength{\tfwidth}{\the\wd\figBox}
\setlength{\tfheight}{\the\ht\figBox}
\setlength{\tfxa}{\tfwidth*\real{#3}}%
\setlength{\tfxb}{\tfwidth*\real{#4}}%
\setlength{\tfya}{\tfheight*\real{#5}}%
\setlength{\tfyb}{\tfheight*\real{#6}}%
\trimFigHeightWithBox{#1}{#2}{\tfxa}{\tfya}{\tfxb}{\tfyb}%
}
\usepackage{xargs}

\newcommandx{\figByHeight}[9][5=0, 6=0, 7=0, 8=0,9=]{
\draw (#1,#2) node[anchor=south west,xshift=-16pt,yshift=-4pt] {\trimh{#3}{#4}{#5}{#6}{#7}{#8}};}
\newcommandx{\figByHeightb}[9][5=0, 6=0, 7=0, 8=0,9=]{
\draw (#1,#2) node[anchor=south west,xshift=-16pt,yshift=-4pt] {\trimhb{#3}{#4}{#5}{#6}{#7}{#8}};}

\newcommandx{\figByHeightWithLabel}[9][5=0, 6=0, 7=0, 8=0,9=]{
\draw (#1,#2) node[anchor=south west,xshift=-16pt,yshift=-4pt] {\trimh{#3}{#4}{#5}{#6}{#7}{#8}} node[draw=white,fill=white,inner sep=1pt,anchor=south west] {#9};}
\newcommandx{\figByHeightWithLabelb}[9][5=0, 6=0, 7=0, 8=0,9=]{
\draw (#1,#2) node[anchor=south west,xshift=-16pt,yshift=-4pt] {\trimhb{#3}{#4}{#5}{#6}{#7}{#8}} node[draw=white,fill=white,inner sep=1pt,anchor=south west] {#9};}

\newcommandx{\figByWidth}[9][5=0, 6=0, 7=0, 8=0,9=]{
\draw (#1,#2) node[anchor=south west,xshift=-16pt,yshift=-4pt] {\trimw{#3}{#4}{#5}{#6}{#7}{#8}};}
\newcommandx{\figByWidthb}[9][5=0, 6=0, 7=0, 8=0,9=]{
\draw (#1,#2) node[anchor=south west,xshift=-16pt,yshift=-4pt] {\trimwb{#3}{#4}{#5}{#6}{#7}{#8}};}

\newcommandx{\figByWidthWithLabel}[9][5=0, 6=0, 7=0, 8=0,9=]{
\draw (#1,#2) node[anchor=south west,xshift=-16pt,yshift=-4pt] {\trimw{#3}{#4}{#5}{#6}{#7}{#8}} node[draw=white,fill=white,inner sep=1pt,anchor=south west] {#9};}
\newcommandx{\figByWidthWithLabelb}[9][5=0, 6=0, 7=0, 8=0,9=]{
\draw (#1,#2) node[anchor=south west,xshift=-16pt,yshift=-4pt] {\trimwb{#3}{#4}{#5}{#6}{#7}{#8}} node[draw=white,fill=white,inner sep=1pt,anchor=south west] {#9};}

\usepackage{todonotes}

\definecolor{jwbGreen}{rgb}{0, .6, 0}
\definecolor{jbaPurple}{HTML}{6600FF}
\definecolor{purple}{rgb}{.7, 0., .8}
\definecolor{pinegreen}{rgb}{0.0, 0.47, 0.44}
\newcommand{\black}{\color{black}}

\usepackage{helvet}

\usepackage{bm}
\usepackage{titlesec}
\usepackage[symbols,nogroupskip,sort=none]{glossaries-extra}
\usepackage{algorithm}
\usepackage{algpseudocode}

\newcommand{\be}{\begin{equation}}
\newcommand{\ee}{\end{equation}}
\newcommand{\bd}{\begin{displaymath}}
\newcommand{\ed}{\end{displaymath}}

\newcommand{\bbb}{{\bm b}}

\newcommand{\pb}{{\bm p}}

\newcommand{\yb}{{\bm y}}

\newcommand{\qb}{{\bm q}}

\newcommand{\wb}{{\bm w}}
\newcommand{\xb}{{\bm x}}
\newcommand{\rb}{{\bm r}}

\newcommand{\whiA}{{\ensuremath{A}}}
\newcommand{\lapL}{\ensuremath{L}}
\newcommand{\matW}{\ensuremath{W}}
\newcommand{\keig}{\ensuremath{k_{\rm eig}}}

\newcommand{\bx}{\mathbf{x}}
\newcommand{\etab}{\boldsymbol{\eta}}

\glsxtrnewsymbol[description={System size}]{N}{\ensuremath{N}}
\glsxtrnewsymbol[description={Number of deflated eigenvalues}]{keig}{\keig}

\glsxtrnewsymbol[description={WaveHoltz system matrix $(N\times N)$}]{whiA}{\whiA}
\glsxtrnewsymbol[description={Laplacian system matrix $(N\times N)$}]{lapL}{\lapL}
\glsxtrnewsymbol[description={Laplacian system matrix $(N\times \keig)$}]{matW}{\matW}

\glsxtrnewsymbol[description={Deflation system matrix $(\keig \times \keig)$}]{WTAW}{$\matW^T\whiA\matW $}

\glsxtrnewsymbol[description={Discretization size in each dimension}]{n}{\ensuremath{n}}
\glsxtrnewsymbol[description={Frequency}]{omega}{\ensuremath{\omega}}
\glsxtrnewsymbol[description={Dimension}]{d}{\ensuremath{d}}

\begin{document}

\begin{frontmatter}
\title{Numerical Study of Eigenvector Deflation to Accelerate the WaveHoltz Method}

\author[vtu]{Daniel Appel\"o\fnref{DanielThanks}}
\ead{appelo@vt.edu}

\address[vtu]{Department of Mathematics, Virginia Tech, Blacksburg, VA 24061 U.S.A.}

\author[rpi]{William D.~Henshaw\corref{cor}\fnref{NSFgrants}}
\ead{henshw@rpi.edu}
\address[rpi]{Department of Mathematical Sciences, Rensselaer Polytechnic Institute, Troy, NY 12180, USA}

\author[hkust]{Zhichao Peng\fnref{ZhichaoThanks}}
\ead{pengzhic@ust.hk}
\address[hkust]{Department of Mathematics, The Hong Kong University of Science and Technology, Hong Kong.}

\cortext[cor]{Corresponding author}

\fntext[NSFgrants]{Research supported by the National Science Foundation under grants DMS-2513122.}

\fntext[DanielThanks]{Research supported by National Science Foundation under grant DMS-2345225, DMS-2436319, and Virginia Tech. This material is based upon work supported by the National Science Foundation under Grant No.
DMS-2424139 while the second author were in residence at the Simons Laufer Mathematical Sciences
Institute in Berkeley, California, during the Fall 2025 semester.}

\fntext[ZhichaoThanks]{Supported in part by the Hong Kong Research Grants Council grants Early Career Scheme 26302724 and General Research Fund 16306825.}

\begin{abstract}
We present a numerical study of eigenvector deflation as a means of accelerating the WaveHoltz method for solving the Helmholtz equation. For energy-conserving (Dirichlet or Neumann) boundary conditions the WaveHoltz fixed-point iteration converges slowly at high frequency, requiring approximately $\mathcal{O}(\omega^{2d})$ iterations in $d$ dimensions. We show that deflating the eigenvectors whose eigenvalues lie nearest the driving frequency substantially reduces iteration counts, and we examine two ways of incorporating the eigenvectors: direct eigenvector deflation (DEVD), in which the forcing and iterate are projected against the deflation set, and augmented-Krylov eigenvector deflation (AUKED) using deflated conjugate gradient (DCG), augmented GMRES (AGMRES), and augmented (recycled) BICGSTAB (ABICGSTAB). The required eigenpairs can be computed efficiently with the EigenWave approach, and we demonstrate, in two dimensions, that when the number of deflation vectors grows quadratically with $\omega$ the asymptotic convergence rate remains essentially constant. Because the eigenvectors on structured grids are naturally represented as matrices, we further apply SVD-based compression to reduce their storage. Numerical experiments on single curvilinear grids discretized with summation-by-parts operators, and on overset grids illustrate the robustness and efficiency of the approach, with the deflated solver breaking even against the undeflated solver after as few as two right-hand sides, when accounting for the cost of precomputing the eigenvectors.
\end{abstract}


\end{frontmatter}

\tableofcontents

\section{Introduction\label{sec:introduction}}

In many engineering applications it is important to solve Helmholtz problems for many right hand sides. If the frequency is high, 
solving the Helmholtz equation becomes challenging due to {the grid} resolution needed to control the dispersion errors \cite{KreOli72} (also called pollution errors \cite{babuska1997pollution}) and the indefinite nature of the linear system of equations arising from discretizing Helmholtz equation \cite{erlangga2008advances,ernst2011difficult,gander2019class}. 

WaveHoltz is an iterative method for solving the Helmholtz equation {based on time filtering the solution to the associated wave equation~\cite{appelo2020waveholtz}}. 
It {results in} a transformed linear system of equations that is {positive definite or close to positive definite, and thus easier to solve by iterative methods than the original indefinite system that arises from the Helmholtz problem}. 
In this work we study how deflation techniques, {using eigenvectors}, (which can also be interpreted as preconditioners) can be used to accelerate the WaveHoltz method. 
Here we consider energy-conserving boundary conditions, such as zero Dirichlet or Neumann boundary conditions. For these the basic WaveHoltz method, {accelerated with a Krylov method such as conjugate-gradients or GMRES is observed to} converge in approximately $O(\omega^d)$ iterations, where $d$ is the space dimension of the problem. 
{Deflation techniques can reduce the number of iterations significantly.}
We note that the we have {previously} used, but not provided a systematic study of, eigenvector deflation techniques in \cite{overHoltzPartOne} where we consider overset grid methods {with implicit time-stepping which results in an order $N$ solver at fixed frequency, where $N$ is the number of grid points.}

A good preconditioner should not only result in low iteration counts but it should also be inexpensive and fast to set up. The number of iterations of deflation preconditioners is reduced when the number of vectors used in the deflation process is increased. However, at the same time, especially when using eigenvectors as we do here, the computation and storage of these can be sizable. Here we employ our EigenWave, \cite{EigenWave}, approach for finding the eigenvectors and show in numerical examples that when finding a large number of eigenvectors the overhead is manageable. In fact we observe that in some cases the ``break-even" of the preconditioner is already at as little as two right hand sides. To reduce the memory footprint when using hundreds of deflation vectors we explore {compression based on the singular value decomposition (SVD)}. As we are using structured grids, {the eigenvectors can be represented as matrices (two dimensional arrays) and use of the SVD} is natural and also turns out to be quite efficient, especially for fine grids. 
{Additionally we evaluate the use of approximate eigenvectors, such as those computed on a coarser grid, or at a lower order of accuracy.}

Deflation for the conjugate gradient (CG)  method was first proposed in \cite{nicolaides1987deflation,dostal1988conjugate} to accelerate its convergence. The method we use here is that of Saad et al. \cite{saad2000deflated}. The augmented GMRES (AGMRES) method appears in different forms and the one we use here is the algorithm given in Baglama and Reichel~\cite{augmentedGMRES2007}. We also use the augmented (recycled) BICGSTAB algorithm (ABICGSTAB) from Amritkar et alia~\cite{Amritkar2015}.  We note that while this is the first systematic investigation of deflation for WaveHoltz, deflation techniques have been {widely used to} accelerate iterative solvers for the Helmholtz equation. For example, multilevel deflation solvers based on complex shifted Laplacian preconditioners (CSLP) are developed in \cite{erlangga2008advances,sheikh2016accelerating}, and further analyzed in \cite{sheikh2013convergence,garciia2018spectrum}. Other aspects of deflation for Helmholtz are considered in \cite{dwarka2020scalable,dwarka2022scalable,chen2023matrix}. 

The rest of the paper is organized as follows. In Section \ref{sec:prelim} we introduce the WaveHoltz method. In Section \ref{sec:eigenvectorDeflation} we review the different deflation algorithms we use. In Section \ref{sec:numerical} we present experiments on single curvilinear meshes, {including results on compression of the eigenvectors, while in} in Section \ref{sec:oversetGridResults} we present numerical results on overset grids. {Concluding remarks are provided in Section \ref{sec:conclusion}}.

\newcommand{\ww}{w}
\newcommand{\vi}{v}

\section{Preliminaries} \label{sec:prelim}
In this section we describe the WaveHoltz method, provide an overview of discrete approximations, and summarize the convergence analysis.

\subsection{The WaveHoltz method for solving the Helmholtz equation}\label{sec:waveholtz}
We are interested in finding numerical approximations to the solution $u=u(\xv)\in\Real$ of a Helmholtz boundary value problem (BVP)
\bse
\label{eq:helmholtzBVP}
\bat
    &  L u  + \omega^2 u = f(x), \qquad&& \bx \in \Omega, \label{eq:helmholtz} \\
    & \Bc u(\bx) = 0,                                           \qquad && \bx \in \partial \Omega.
\eat
\ese
The operator $L$ in~\eqref{eq:helmholtzBVP} is taken here as
\ba
   L \eqdef \nabla \cdot (c(\bx)^2 \nabla).   \label{eq:L}
\ea
The function $f(\xv)\in\Real$ is a given forcing function, $c(\xv)>0$ is a given variable coefficient, and $\omega \ge 0$ is the given frequency. 
The boundary condition operator $\Bc$ is taken here to be a Dirichlet condition although similar results will hold for a Neumann boundary condition.
The BVP~\eqref{eq:helmholtzBVP} will have a solution if the problem is not at resonance, that is 
provided $\omega$ is not an eigenvalue of the corresponding eigenvalue problem.
The WaveHoltz method can also be applied to problems with radiation boundary conditions or damping~\cite{Rotem2024,WaveHoltz2} where the solution $u$ would be complex valued, but
deflation for problems of this type are not considered here.

The WaveHoltz method~\cite{appelo2020waveholtz} is an iterative solver to find solutions to~\eqref{eq:helmholtzBVP} 
that is built upon a time-domain solver for the associated wave equation initial boundary value problem (IBVP)
\begin{subequations}
\label{eq:waveEquation}
\begin{align}
    &\ww_{tt}= L \ww -f(\bx)\cos(\omega t),\\
    &\ww(\bx,t)=0,\ \ \bx\in \partial\Omega,\\
    &\ww(\bx,0)=\vi^{(k)}(\bx), \label{eq:waveEquationIC}\\     
    &\ww_t(\bx,0)=0.
\end{align}
\label{eq:wave}
\end{subequations}

The WaveHoltz algorithm iterates on the initial condition $\vi^{(k)}$ in~\eqref{eq:waveEquationIC} to find a time-periodic solution. Given a current guess, $\vi^{(k)}$, the next iterate, $\vi^{(k+1)}$ is found by solving the wave equation for $w(\xv,t;\vi^{(k)})$ and then 
time filtering this solution, 
\begin{align}
   \vi^{(k+1)} =  \Pi \,\vi^{(k)} \eqdef \frac{2}{\Tbar}\int_{0}^{\Tbar} \big( \cos(\omega t)-\frac{1}{4} \big) \, \ww(\bx,t; \vi^{(k)}) \, dt
     {~= S v^{(k)} + b,} 
     \label{eq:waveHoltzFilter}
\end{align}
where $\Tbar= \Np T$, $T=2\pi/\omega$ is the period, with $\Np$ is a positive integer. The affine operator $\Pi$ filters the time-domain solution 
 over $N_p$ periods {and can be written as $\Pi \, v^{(k)} = S v^{(k)} + b$ where $S$ is a linear operator and $b$ depends in the forcing $f$}.
As shown in \cite{appelo2020waveholtz}, the solution $u$ to the Helmholtz BVP~\eqref{eq:helmholtzBVP} is a fixed point of the filtering operator $\Pi$ and the corresponding time-domain solution is $\ww(\bx,t)=u(\bx)\cos(\omega t)$. In other words, solving the Helmholtz problem~\eqref{eq:helmholtzBVP}
is equivalent to solving the fixed point problem
\begin{equation}
    \Pi u = u.
\end{equation}
This fixed point problem can be further rewritten as an equivalent linear problem 
\begin{align}
    \Ac u = (I-S)u = b  ,\label{eq:I_minus_S}
\end{align}
where $\Ac\eqdef I-S$, and $b \eqdef\Pi \, 0$ denotes the result of solving the IBVP~\eqref{eq:waveEquation} with zero initial conditions.
The linear operator $S$ is thus defined as $Su = \Pi\, u -\Pi\, 0$.

\renewcommand{\algFontSize}{\small}
\begin{algorithm}[t]
\algFontSize 
\caption{WaveHoltz Algorithm - Fixed-Point Iteration.}
\begin{algorithmic}[1]

  \Function{WaveHoltz}{$\omega$,$f$,$\Np$}  
    \State // Final time is $\Tbar=\Np T$, with $T=2\pi/\omega$.  \label{alg:init}
    \State $k=0$ \Comment WaveHoltz~iteration counter.
    \State $v^{(k)}=0$   \Comment Assign initial guess for Helmholtz iterate 
       \label{alg:v0}
    \While{ not converged} \Comment Start WaveHoltz~iterations.

      \State $w^{(k)}(\xv,0) = v^{(k)}(\xv)$ \Comment Initial condition for wave equation solve.
      \State $w^{(k)}(\xv,0:\Tf)$ = \Call{solveWaveEquation}{$w^{(k)}(\xv,0)$,$f$} \Comment Solve~\eqref{eq:waveEquation} for $\wv(\xv,t)$, $t\in[0,\Tf]$.\label{alg:solveWave}
      \State $\displaystyle v^{(k+1)}(\xv) = \f{2}{\Tbar} \int_{0}^{\Tbar} \left( \cos(\omega t) - \f{1}{4} \right) \, w^{(k)}(\xv,t; v^{(k)}) \, dt$
           \Comment Time filter the wave equation solution.
      \State $k = k+1$
    \EndWhile    \Comment End WaveHoltz iterations.
    \State $u(\xv) = v^{(k)}(\xv)$ \Comment Approximate Helmholtz solution.
 \EndFunction
\end{algorithmic} 
\label{alg:waveHoltz}
\end{algorithm}

Algorithm~\ref{alg:waveHoltz} gives the WaveHoltz fixed-point iteration. 
This algorithm can be accelerated using traditional methods such as Chebyshev acceleration or modern Krylov sub-space methods. 
{Note that when} solved as a linear system of equations (\ref{eq:I_minus_S}) the application of $(I-S)$ by Algorithm~\ref{alg:waveHoltz} can be carried out with $f=0$. 

Let us now introduce the idea of accelerating the WaveHoltz method using eigenfunctions (eigenvectors in the discrete case), which is called eigenvector deflation. 
The eigenvalue problem corresponding to~\eqref{eq:helmholtzBVP} is
\bse
\label{eq:eigenvalueBVP}
\bat
    &  L \phi = - \lambda^2 \phi  \qquad&& \bx \in \Omega, \label{eq:eigenvaluePDE} \\
    & \Bc \phi(\bx) = 0,                                           \qquad && \bx \in \partial \Omega.
\eat
\ese
Assume that~\eqref{eq:eigenvalueBVP} has discrete eigenvalues $\lambda_m>0$ 
and a complete set of orthonormal eigenfunctions $\phi_m$, $m=1,2,\ldots$.


\renewcommand{\algFontSize}{\small}
\begin{algorithm}[b]
\algFontSize 
\caption{WaveHoltz Algorithm with Eigenfunction Deflation.\label{alg:waveHoltzDirectDeflation}}
\begin{algorithmic}[1]

  \Function{WaveHoltz}{$\omega$,$f$,$\Np$}  
    \State // Final time is $\Tbar=\Np T$ where $T=2\pi/\omega$. 
    \State $k=0$ \Comment WaveHoltz~iteration counter.
    \State $v^{(k)}=0$   \Comment Initial guess for Helmholtz iterate (deflate if non-zero) 
    \State \blue $ \displaystyle  f_d =  f - \sum_{\phi_m\in\DeflationSet} ( f, \phi_m )_\Omega \, \phi_m$ \Comment Deflate forcing. \black \label{eq:deflateForcing}
    \While{ not converged} \Comment Start WaveHoltz~iterations.

      \State $w^{(k)}(\xv,0) = v^{(k)}(\xv)$ \Comment Initial condition for wave equation solve.
      \State $w^{(k)}(\xv,0:\Tf)$= \Call{solveWaveEquation}{$w^{(k)}(\xv,0)$,$f_d$} \Comment Solve for $\wv(\xv,t)$ using deflated forcing $f_d$. 
      \State $\displaystyle v^{(k+1)}(\xv) = \f{2}{\Tbar} \int_{0}^{\Tbar} \left( \cos(\omega t) - \f{\alpha}{2} \right) \, w^{(k)}(\xv,t) \, dt$
           \Comment Time filter.
      \State \blue $ \displaystyle v^{(k+1)} =  v^{(k+1)} - \sum_{\phi_m\in\DeflationSet} ( v^{(k+1)}, \phi_m )_\Omega \, \phi_m$ \Comment Deflate iterate (skip with true eigenfunctions). \black
         \label{eq:deflationIterate}
      \State $k = k+1$
    \EndWhile    \Comment End WaveHoltz iterations.
    \State \blue $\displaystyle v^{(k)} =  v^{(k)} + \sum_{\phi_m\in\DeflationSet} \f{(f,\phi_m)_\Omega}{\omega^2-\lambda_m^2} \,  \phi_m  $ \Comment Inflate. \black
        \label{eq:inflateSolution}
    \State $u = v^{(k)}$; \Comment Approximate Helmholtz solution.
 \EndFunction
\end{algorithmic} 
\end{algorithm}

The WaveHoltz algorithm in Algorithm~\ref{alg:waveHoltz} can be accelerated using some subset $\phi_m \in \DeflationSet$
of the eigenfunctions, where $\DeflationSet$ is called the \textsl{deflation set}. We will denote the dimension of this set by $N_{\DeflationSet}$.
The choice of which eigenfunctions are best to use is discussed in Section~\ref{sec:discreteApproximations}.
The WaveHoltz method with eigenfunction deflation is given in Algorithm~\ref{alg:waveHoltzDirectDeflation}.
There are two key changes to Algorithm~\ref{alg:waveHoltz}.
The first is on line~\ref{eq:deflateForcing} where the forcing $f$ is adjusted by removing the components of $f$ along the eigenfunctions in the deflation set.
Here $(f,\phi_m)_\Omega$ denotes the inner product of $f$ and $\phi_m$.
The second key change is on line~\ref{eq:inflateSolution} where the components of the solution in the deflation set are added to the WaveHoltz solution.
Line~\ref{eq:deflationIterate} contains an additional change where the WaveHoltz iterate $v^{(k+1)}$ is deflated after each iteration. 
This line is useful to add when approximate eigenvectors are being used. For example, 
on a overset grid the discrete eigenvectors may not be orthogonal
to machine precision in which case using Line~\ref{eq:deflationIterate} is useful.
There is an alternative way to incorporate approximate eigenvectors into a deflation process and this uses augmented Krylov methods;
this is discussed further in Section~\ref{sec:eigenvectorDeflation}.

\subsection{Spatial and temporal discretization}

We have presented the methodology of the WaveHoltz method at the continuous level.
In fact, it can be flexibly integrated with various time-domain solvers for the wave equation such as those based on finite-difference, finite-volume and finite element approximations. To provide a concrete example, consider using a finite-difference approximation on a structured curvilinear grid.
Let $\xv_\jv$ denote the grid-points, where $\jv=[j_1,j_2,j_3]$ is a multi-index with $j_m=0,1,2,\ldots,N_m$, 
and $N_m$ is the number of grid cells in coordinate direction $m$.
Let $U_\jv \approx u(\xv_\jv)$ denote the discrete approximation to the Helmholtz solution.
A $p$-th order accurate approximation to the Helmholtz BVP~\eqref{eq:helmholtzBVP} can be written as
\bse
 \label{eq:discreteHelmholtz}
\bat
    & \Lph U_\jv + \omega^2 U_\jv = f(\xv_\jv)   , \qquad && \jv\in\Omega_h, \\
    & \Bv_h U_\jv = 0                                \qquad && \jv\in\p\Omega_h,  \label{eq:discreteHelmholtzBCs}
\eat
\ese
where $\Omega_h$ denotes the set of grid points where the interior equation is applied and $\p\Omega_h$ denotes the
set of boundary grid points. 
{$L_{p,h}$ denotes a a $p$-th order accurate difference operator.}
The discrete boundary equations in~\eqref{eq:discreteHelmholtzBCs} represent one or more discrete boundary
conditions as needed for a $p$-th order accurate approximation, for example the Dirichlet condition could be augmented by compatibility conditions~\cite{lcbc2022,wimp2025}.
{
The corresponding discrete eigenvalue problem is 
\bse
 \label{eq:discreteEigenvalueProblem}
\bat
    & \Lph U_\jv = - \lambda_h^2 \Phi_\jv   , \qquad && \jv\in\Omega_h, \\
    & \Bv_h \Phi_\jv = 0                    \qquad && \jv\in\p\Omega_h.
\eat
\ese
Let us suppose that~\eqref{eq:discreteEigenvalueProblem} has positive eigenvalues $\lambda_{h,m}>0$, and
a complete set of orthonormal eigenvectors $\Phiv_{m,h}$, for $m=1,2,\ldots,N$, where $N$ denotes the number of degrees of freedom in $\Omega_h$.
Eigenvalues and eigenvectors of~\eqref{eq:discreteEigenvalueProblem} will be used for deflation.
}

The WaveHoltz method will be used to find the solution $U_\jv$ to the discrete BVP~\eqref{eq:discreteHelmholtz}.
To do this, a discrete approximation to the IBVP for the wave equation~\eqref{eq:waveEquation} will be constructed. 
For example, an explicit finite-difference discretization for the wave equation, second-order accurate in time and $p$-th order accurate in space, is 
\bse
\label{eq:discreteWaveExplicit}
\ba
   & \Dpt\Dmt W_\jv^n = \Lph W^n_\jv - f(\xv_\jv) \, \cos(\omega t^n),   \qquad&& \jv\in\Omega_h, \quad n=0,1,2,\ldots,N_t, \label{eq:discreteWaveExplicitInterior}\\
   & W_\jv^0 = V_\jv^{k} ,                                                 \qquad&& \jv\in\Omega_h, \\
   & \Dzt W_\jv^0 = 0 ,                                                    \qquad&& \jv\in\Omega_h, \\
   & \Bv_h W_\jv^n = 0  ,                                                  \qquad && \jv\in\p\Omega_h, 
\ea
\ese
where $W_\jv^n \approx w(\xv_\jv,t^n)$, 
$t^n = n \dt$, $\dt$ is the time-step,
and $N_t$ is the number of time-steps. 
Here $\Dzt, \Dpt$ and $\Dmt$ are the standard centered, forward, and backward difference operators in time, 
with 
$\Dzt W_\jv^n =(W_\jv^{n+1}- W_\jv^{n-1})/(2\dt)$,  $\Dpt W_\jv^n =(W_\jv^{n+1}- W_\jv^n)/\dt$, and $\Dmt W_\jv^n =(W_\jv^{n}- W_\jv^{n-1})/\dt$.
On the other hand, an unconditionally stable implicit time-stepping method replaces~\eqref{eq:discreteWaveExplicitInterior} with 
\ba
   \Dpt\Dmt W_\jv^n = \half \Lph\Big[  W^{n+1}_\jv + W^{n-1}_\jv \Big]  - f(\xv_\jv) \, \cos(\omega t^n). \label{eq:discreteWaveImplicit}
\ea
Implicit methods can be useful since a very large time-step can be chosen, with as few as five time-steps per-period, and this can lead
to an $O(N)$ WaveHoltz algorithm at fixed frequency, where $N$ denotes the number of grid-points~\cite{overHoltzPartOne}.
The time filter~\eqref{eq:waveHoltzFilter} can be approximated with the trapezoidal rule,
\begin{align}
  \Pi_h V_\jv^{(k)} = \frac{2}{\Tbar} \sum_{n=0}^{N_t}\eta_n \big(\cos(\omega t)-\frac{1}{4}\big) W_\jv^n \dt,\quad 
    \eta_n=\begin{cases}
            \half,\quad n=0\;\text{or}\; n=N_t,\\
            1,\quad\text{otherwise}, 
            \end{cases}
            \label{eq:quadrature}
\end{align}
which is spectrally accurate for periodic functions.
In matrix-vector form the discrete WaveHoltz fixed-point iteration is then
\bse
\label{eq:WaveHoltzDiscreteFPI}
\ba
  & \Vv^{(k+1)} = \Pi_h \Vv^{(k)} = S_h \Vv^{(k)} + \bv,  \label{eq:Sh} \\
  & \bv \eqdef  \Pi_h \, \zerov,
\ea
and the discrete Helmholtz solution satisfies
\ba
   A \Uv \eqdef (I - S_h) \Uv = \bv, \label{eq:WaveHoltzMatrixEquation}
\ea
\ese
where $\Vv^{(k)}$ and $\Uv$ denote vectors of all the grid point values $V_\jv^{(k)}$ and $U_\jv$, respectively.
Note that the matrix $S_h$ need not be explicitly formed, rather its application on a vector is computed by time-stepping the wave equation
and filtering the result (called a wave-solve).

The approximations in~\eqref{eq:discreteWaveExplicit} and~\eqref{eq:discreteWaveImplicit} are second-order accurate in time.
The WaveHoltz solutions computed using these approximations would have $O(\dt^2)$ time discretization errors. 
However, these time-discretization errors can be completely removed, by for example, using an adjusted value for $\omega$ in~\eqref{eq:discreteWaveExplicit}.
There are multiple ways to remove the time-discretization errors and details of these approaches can be found in~\cite{stolk2021time, peng2022waveholtz, ElWaveHoltz, overHoltzPartOne}. In this paper we use the approach from \cite{overHoltzPartOne}.

Pollution (dispersion) errors are a serious problem for high frequencies and require the use of a much finer grid than might be expected.
An analysis in~\cite{overHoltzPartTwo} gives a rule-of-thumb for estimating the number of points per wave-length (PPW) for 
a $p$-th order accurate finite difference approximation ($p$ an even positive integer).
{
Given $\omega$, the domain size $\domainSize$, and a relative error tolerance $\eps$, 
the PPW rule of thumb gives the number of points per wavelength for a $p^{\rm th}$-order accurate scheme ($p$ an even integer) to be 
\bse
\ba
  \PPW_p \eqdef  2 \pi \, (\pi \, b_{p/2})^{1/p} \, \left[ \f{N_\Lambda}{\eps} \right]^{1/p}, 
  ~~ N_\Lambda\eqdef \f{\domainSize}{\Lambda}, 
  ~~\Lambda \eqdef \f{2\pi c}{\omega},
  ~~b_{\mu}  \eqdef \f{2\, (\mu !)^2}{ (2\mu+2)!}, 
   \label{eq:ppwPollution}
\ea
where $N_\Lambda$ is the size of the domain in wavelengths.
For a Cartesian grid, this implies that the number of grid points, in a given direction, should be chosen as 
\ba
   N_x = \f{\domainSize}{\dx} = \omegaN \, \omegaN^{1/p} \left[ \f{b_{p/2}}{2 \eps} \right]^{1/p} , \label{eq:NxFromPPW}
\ea
where $\omegaN$ denotes the non-dimensional frequency
\ba
  &  \omegaN \eqdef \f{\omega \domainSize}{c} .   \label{eq:omegaN}
\ea
\ese}

\newcommand{\betad}{\beta_d}
\subsection{Convergence theory for the WaveHoltz method} \label{sec:discreteApproximations}


%
The convergence theory for the continuous and discrete WaveHoltz algorithm is given in~\cite{appelo2020waveholtz}, see also~\cite{overHoltzPartOne}, for further discrete analysis. 
Here we summarize results of the theory that will be germane to understanding the convergence of the WaveHoltz iterations using 
deflation and Krylov methods.
The convergence theory is based on eigenfunction (eigenvector) expansions.
In terms of the eigenvalues in~\eqref{eq:discreteEigenvalueProblem}, the eigenvalues of the WaveHoltz operator $S_h$ in~\eqref{eq:Sh} are 
\bse
\ba
    \beta_m = \beta(\lambda_{h,m};\omega), 
\ea
where the WaveHoltz filter function $\beta(\lambda;\omega)$ is 
\ba
  \beta(\lambda;\omega) & \eqdef \frac{2}{\Tbar}\int_{0}^{\Tbar}(\cos(\omega t)-\frac{1}{4})\cos(\lambda t) \, dt  , \\
                        & = \sinc(\omega+\lambda,\Tbar) + \sinc(\omega-\lambda,\Tbar) - \half \, \sinc(\lambda,\Tbar) .  \label{eq:beta_fun}
\ea
\ese
The asymptotic convergence rate (ACR) of the WaveHoltz FPI~\eqref{eq:WaveHoltzDiscreteFPI} is then
\ba
   \rho_{\ssf\rm FPI} = \max_{m} | \beta(\lambda_{h,m}; \omega) | , \label{WaveHoltzFPIconvergence}
\ea
where we have ignored various corrections due to the
finite time-step and finite grid spacings~\cite{overHoltzPartOne}. We note that $\beta(\lambda;\omega)$ is only a function of $\lambda/\omega$ when $T=2\pi/\omega$.
The left graph in Figure~\ref{fig:waveHoltzBeta} plots $\beta$ for different values of $\Np$. 
The convergence rate of the FPI will generally be determined by $\beta$ evaluated at the closest eigenvalue $\lambda_{h,m}$ to $\omega$.
When using deflation, the ACR will be determined by~\eqref{WaveHoltzFPIconvergence} but with the maximum taken over the set of eigenvalues excluding
those in the deflation set, $\DeflationSet$. Therefore one should normally deflate the eigenvectors with eigenvalues closest to $\omega$.

As discussed further in Section~\ref{sec:eigenvectorDeflation}, the convergence rate of Krylov methods such as conjugate gradient and GMRES will depend on the condition number, $\kappaA$, of the matrix $A$ in~\eqref{eq:WaveHoltzMatrixEquation}. Since $A$ is symmetric, or close to symmetric, for the discretizations we consider, this condition number can be approximated by  
\bse
\ba
  \kappaA = \f{\max_m{|\mu_{h,m}|}}{\min_m{|\mu_{h,m}|}}, 
\ea
where $\mu_{h,m}$ are the eigenvalues of $A$
\ba
   \mu_{h,m} \eqdef \mu(\lambda_{h,m}; \omega) = 1  - \beta(\lambda_{h,m}; \omega).  \label{eq:eigenvaluesOfM}
\ea
\ese
The right graph of Figure~\ref{fig:waveHoltzBeta} plots the function $\mu(\lambda;\omega) = 1-\beta(\lambda;\omega)$.
As with the FPI, for Krylov methods one should normally deflate the eigenvectors with eigenvalues closest to $\omega$, and these correspond
to the smallest eigenvalues of $A$.

{
\newcommand{\figw}{10cm}
\newcommand{\figh}{3.8cm}
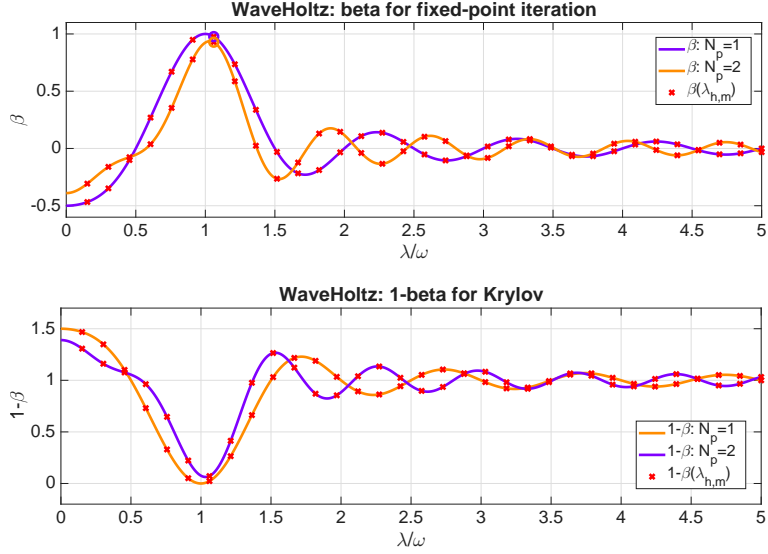
\begin{figure}[htb]
\begin{center}
\begin{tikzpicture}
  \useasboundingbox (0,.5) rectangle (\figw,2*\figh);  
  \begin{scope}[yshift=0*\figh]
    \figByWidth{0.0}{1*\figh}{fig/waveHoltzBetaFunctionPlusEigs}{\figw}[0.][0.][0.][0.]
    \figByWidth{0.0}{0*\figh}{fig/waveHoltzOneMinusBetaFunction}{\figw}[0.][0.][0.][0.]
  \end{scope}  
\end{tikzpicture}
\end{center}
\caption{
  Top: WaveHoltz filter function $\beta$ for $\Np=1$, and $\Np=2$ periods per time-interval.
  The asymptotic convergence rate of the fixed-point iteration corresponds to the circled value.
  Bottom: Plots of $1-\beta(\lambda)$ for $\Np=1$, and $\Np=2$ periods per time-interval.
   The Krylov solvers operate on the matrix $A=I-S_h$ which has eigenvalues $1-\beta(\lambda_{h,m})$ of which representative values are shown with red x's.
   }
\label{fig:waveHoltzBeta}
\end{figure}
}

Of particular interest is the convergence behaviour of WaveHoltz for large $\omega$.
As noted above, the WaveHoltz convergence rate is normally determined by the nearest eigenvalues to $\omega$, 
and it can be shown that for $\lambda$ near $\omega$, 
\bse
\ba
  & \beta(\lambda;\omega) = 1 - C_\beta \, z^2 + O(z^3), \qquad  C_\beta \eqdef \f{2}{3} N_p^2 \, \pi^2 -\f{1}{4}, \label{eq:betaGap}
\ea
where $z$ is the \textsl{gap width},
\ba
  z\eqdef \f{\lambda-\omega}{\omega}. 
\ea
\ese
Using~\eqref{eq:betaGap} in~\eqref{WaveHoltzFPIconvergence} leads to an estimate of the ACR for the FPI
\bse
\ba  
    & \rho_{\ssf\rm FPI} \approx 1 - C_\beta \, \delta^2 , \qquad
     \delta \eqdef \min_m \left| \f{\lambda_{h,m}-\omega}{\omega} \right|. 
\ea 
\ese
Weyl~\cite{weyl1912asymptotische} studied the distribution of eigenvalues of the Laplacian, $L=\Delta$, with Dirichlet boundary conditions, 
and showed that for fairly general domains the gap width scales as $\delta = O(\omega^d)$ in $d$ dimensions for $\omega\rightarrow \infty$.
In this case the ACR for the FPI scales as $\rho_{\ssf\rm FPI} \approx 1 - O(\omega^{2d})$, which implies $O(\omega^{2d})$ iterations for convergence.
On the other hand, using~\eqref{eq:eigenvaluesOfM}, and noting that $|\mu_{h,m}|<1.5$,  the condition number of $A$ will behave as $\kappaA = O(\omega^{2d})$.
From standard convergence estimates for symmetric problems, Krylov methods such as conjugate gradient or GMRES will generally require $O(\omega^{d})$ iterations to
converge. Thus, for increasing $\omega$, a very large number of Krylov, and an even larger number of FPI iterations will generally be needed for convergence. Eigenvector deflation, as discussed in Section~\ref{sec:eigenvectorDeflation}, is one way to decrease the number of iterations.


\section{Eigenvector deflation and deflated Krylov methods}\label{sec:eigenvectorDeflation}

The basic WaveHoltz fixed-point iteration can be effectively accelerated with Krylov sub-space methods.
When using energy conserving Dirichlet or Neumann boundary conditions, however, the convergence can be slow for large $\omega$ since generically
the problem is almost always close to resonance. 
As discussed in Section~\ref{sec:discreteApproximations}, this implies Krylov algorithms such as conjugate-gradient or GMRES will
generally require $O(\omega^d)$ iterations in $d$ dimensions leading to a high computational cost, especially in $d=3$ space dimensions.

Eigenvector deflation can be used reduce the computational cost
by using discrete approximations to some eigenvectors to accelerate the convergence.
This approach consists of two stages. 
In the first stage some number of eigenvectors or approximate eigenvectors are computed.
Approximate eigenvectors could be eigenvectors computed on a coarser grid or using a lower-order accurate approximation.
In the second stage these pre-computed eigenvectors are used within the WaveHoltz algorithm.
Eigenvectors can be incorporated into the WaveHoltz algorithm in multiple ways.
We will consider two ways called \textsl{direct-eigenvector-deflation} (DEVD) and \textsl{augmented-Krylov-eigenvector-deflation} (AUKED).

DEVD makes use of accurate discrete eigenvectors of the eigenvalue problem~\eqref{eq:discreteEigenvalueProblem}.
The continuous analogue of DEVD was given in Algorithm~\ref{alg:waveHoltzDirectDeflation}.
The components of the forcing along the set of eigenvectors in the deflation set are explicitly removed at the start of the WaveHoltz iteration.
The WaveHoltz algorithm is then used to compute a solution to the problem with deflated forcing using a standard FPI or Krylov algorithm. 
Finally the computed solution is corrected by adding in the components of the solution along the eigenvectors in the deflation set.
The WaveHoltz solution can be optionally deflated after each iteration; this can be useful when using approximate eigenvectors that are
accurate to the size of truncation errors. For example, on a overset grid, the discrete eigenfunctions may only be 
orthonormal to the discretization error of a $p$-th order scheme, $O(h^p)$, and in this case deflating the WaveHoltz solution after each iteration is beneficial.

In AUKED, accurate discrete or approximate eigenvectors are used to augment the Krylov sub-space algorithm.
In fact the augmented vectors can be any set of linearly independent vectors.
There are a number of augmented Krylov methods including deflated conjugate-gradient (DCG), 
augmented GMRES (AGMRES), and augmented (recycled) bi-CG-Stab (ABICGSTAB).
DCG can be used for WaveHoltz problems with a symmetric matrix $S_h$.
For overset grid problems where $S_h$ is not symmetric, AGMRES and ABICGSTAB are used.

\subsection{Deflated conjugate gradient (DCG)\label{sec:dcg}}

To accelerate the convergence of the WaveHoltz method for symmetric discretizations, we apply the deflated conjugate gradient (DCG) method \cite{saad2000deflated} by adding eigenvectors associated with the smallest eigenvalues {of the matrix $A$ in~\eqref{eq:WaveHoltzMatrixEquation}} to the Krylov subspace.
{These eigenvalues of $A$ will generally correspond to eigenvalues of the discrete eigenvalue problem~\eqref{eq:discreteEigenvalueProblem} that are close to $\omega$.
}
 
Consider a symmetric positive definite (SPD) matrix $A$ and linear system $A\xb=\bbb\in\mathbb{R}^N$. Given an initial guess $\xb_0$, in the $k$-th iteration, the CG method seeks the solution $\xb_k$ in the Krylov subspace 
\begin{equation}
    \xb_0+ \Kc(\rb_0,\whiA,k)=\xb_0+\{\rb_0,\whiA\rb_0,\dots,\whiA^{k-1} \rb_0\},\quad \rb_0=\bbb-\whiA\xb_0,
\end{equation}
corresponding to the optimal error with respect to the $||\cdot||_A$ norm,  $||\yb||_A=\sqrt{\yb^T\whiA\yb}$. In other words,
\begin{align*}
\xb_k = \min_{\hat{\xb} \in \xb_0+ \Kc(\rb_0,\whiA,k)}||\xb-\hat{\xb}||_A.
\end{align*}
The DCG method augments the Krylov subspace  with a deflation space, which is the column space of a deflation matrix $W\in\mathbb{R}^{N\times \Ndeflate}$, {where $\Ndeflate$ denotes the number of vectors in the deflation space (often chosen to be eigenvectors)}. The DCG method initializes $\xb_0$ in the deflation space such that $$W^TA(\xb_0-\bbb)=W^T\rb_0=0$$ and performs a CG search in the subspace $A$-orthogonal to the deflation space. 
This strategy results in the DCG algorithm \cite{saad2000deflated}, which is presented in Algorithm \ref{alg:DCG}.

As discussed in \cite{saad2000deflated}, the DCG method can be seen as a preconditioned CG method with the preconditioner $PP^T$, where $P$ is the projector onto the $A$-orthogonal complement of the deflation space determined by $W$:
\begin{equation}
    P = I-W(W^TA W)^{-1}W^TA.
\end{equation}
Utilizing this point of view, the convergence speed of the DCG method is analyzed in \cite{saad2000deflated} and summarized in the following theorem. 
\begin{thm}\label{thm:dcg}
(Theorem 4.3 of \cite{saad2000deflated})
Let $\kappaD$ be the condition number of $P^TAP$. The approximate solution of $A\xb=\bbb$ generated by the $j$-th iteration of the DCG algorithm in Algorithm \ref{alg:DCG} satisfies
\begin{equation}
    ||\xb-\xb_j||_A\leq2\left(\frac{\sqrt{\kappaD}-1}{\sqrt{\kappaD}+1}\right)^j|| \xb-\xb_0||_A, \quad ||\yb||_A=\sqrt{\yb^TA\yb}.
\end{equation} 
\end{thm}
\begin{algorithm}
\caption{DCG($\whiA,\bbb,W$): deflated conjugate gradient to solve $\whiA\xb=\bbb$ with the deflation matrix $W$.}\label{alg:DCG}
\begin{algorithmic}[1]
\State Input: $\whiA$, $\bbb$ and $W$.
\State $\xb_{-1}=0$
\State $\rb_{-1} = \bbb-\whiA \xb_{-1}$
\State $\xb_0 =  \xb_{-1} + W (\matW^T \whiA \matW)^{-1} ( \matW^T \rb_{-1} )$
\State $ \rb_0 = \bbb - \whiA \xb_0 $
\State $\etab_0 = \matW^T \whiA \rb_0$
\State $ \muv_0 = (\matW^T \whiA \matW)^{-1} \etab_0$
\State $ \pb_0 = \rb_0 - \matW {\bm \mu_0}$
\State $ r_2 = \rb^T \rb$
\State $ r_0 = r_2$
\State $k = 0$
\While{$r_2/r_0 > {\rm (TOL)}^2$ }  
 \State  $\qb_k  = \whiA \pb_k$
 \State  $\alpha_k \gets r_2/ (\pb_k^T \qb_k)$
 \State  $\xb_{k+1} = \xb_k + \alpha_k \pb_k$
 \State  $ \rb_{k+1} = \rb - \alpha_k \qb$
 \State  $ \hat{r}_2 \gets r_2$
 \State  $ r_2 \gets \rb_{k+1}^T \rb_{k+1}$
  \State  $ \beta_{k+1} \gets  r_2/\hat{r}_2$
  \State {\blue $\etab_{k+1} = \matW^T \whiA \rb_{k+1}$}            \Comment {\blue Skip with true eigenvectors.} \black
  \State  ${\blue \muv_{k+1}= (\matW^T \whiA \matW)^{-1} \etab_{k+1}}$  \black \Comment {\blue Skip with true eigenvectors.} \black
  \State $ \pb_{k+1} = \beta_{k+1} \pb_k+\rb_{k+1}  {\blue - \matW \muv_{k+1} }$ \Comment {\blue Ignore $ -\matW \mu_{k+1}$ with true eigenvectors.}\black
  \State $k \gets k+1$
\EndWhile
\end{algorithmic}
\end{algorithm}

The main modifications to the CG algorithm made by the deflation method are as follows.
Line $3$ of Algorithm \ref{alg:DCG} constructs  an initial guess $\xb_0\in W$ such that  
\begin{align*}
W^T\rb_0 &= W^T(\bbb-\whiA\xb)=W^T\bbb - W^T\whiA\left(\xb_{-1}+AW(W^TAW)^{-1}W^T\rb_{-1}\right)\notag\\
         &=W^T\bbb - W^T\whiA\xb_{-1}+W^TAW(W^TAW)^{-1}W^T\rb_{-1}= W^T\rb_{-1}-W^T\rb_{-1}=0.
\end{align*}
Line $5,6,20$ and $21$ of Algorithm \ref{alg:DCG} enforce the search direction $\pb_{k}$ to be $A$-orthogonal to the deflation space $W$.

When eigenvectors are used as the deflation matrix, after initialization through the projection onto the deflation space which is a span of eigenvectors, the Krylov vectors generated by the CG algorithm is always $A$-orthogonal to each other \cite{saad2000deflated}. In other words, when eigenvectors are chosen as the deflation matrix, line 20-22 in Alg. \ref{alg:DCG} enforcing $A$-orthogonality after each iteration is not necessary. 

The condition number of the eigenvector-deflated system $P^TAP$ is discussed in Sec. 5 of \cite{saad2000deflated}, and we summarize it as a theorem here.
\begin{thm}(Sec. 5 of \cite{saad2000deflated})\label{lem:kappa}
Let the eigenvalues of $A$ be $0< \tilde{\lambda}_1\dots\leq\tilde{\lambda}_N$. If the columns of the deflation matrix 
$$W=\left[\wb_1,\wb_2,\dots,\wb_{N_{\DeflationSet}}\right]$$ are eigenvectors corresponding to $\tilde{\lambda}_1,\dots,\tilde{\lambda}_{N_{\DeflationSet}}$, then the condition number of  $P^TAP$
is 
\ba
   \kappaD \eqdef \frac{\tilde{\lambda}_N}{\tilde{\lambda}_{N_{\DeflationSet}+1}}.
\ea
\end{thm}

{
From the equation for the $\beta$ function in~\eqref{eq:beta_fun}, $\beta\ge -0.5$ and thus 
the largest eigenvalue of the WaveHoltz matrix $A = I - S_h$ is no greater than $1.5$.
Referring to the bottom graph of Figure~\ref{fig:waveHoltzBeta}, suppose all eigenvectors with eigenvalues $1-\beta(\lambda_{h,m}) \le 0.5$ are deflated.
This would then lead to a condition number $\kappaD = 1.5/0.5 = 3$ that is independent of $\omega$ and a convergence rate  
\[
   \frac{\sqrt{\kappaD}-1}{\sqrt{\kappaD}+1} \approx 0.27.
\]
With this convergence rate the residual can be expected to be reduced nine orders of magnitude in 20 iterations or fourteen orders of magnitude in 30 iterations. Engineering precision can likely be reached in $5-7$ iterations.
}

\newcommand{\xvp}{\xv'}
\subsection{Augmented GMRES (AGMRES)} \label{sec:augmentedGMRES}

The augmented GMRES (AGMRES) method is given in Algorithm~\ref{alg:augmentedGmresWithGivens}.
This is the algorithm given in Baglama and Reichel~\cite{augmentedGMRES2007} except that the version in Algorithm~\ref{alg:augmentedGmresWithGivens} makes 
use of Givens rotations to efficiently solve the least squares problem. As with the deflated conjugate gradient method, AGMRES uses augmented vectors $\wv_j$, $j=,1,2,\ldots,N_{\DeflationSet}$,
which may be eigenvectors or approximate eigenvectors.
The changes from the standard GMRES algorithm are highlighted in {\blue blue}.

\renewcommand{\algFontSize}{\small}
\begin{algorithm}
\algFontSize 
\caption{Augmented GMRES algorithm with Givens rotations. }
\begin{algorithmic}[1]

  \Function{$\xv_j=$augmentedGMRES}{$A$, $\bv$, $\xv_0$, $W$, $N_{\DeflationSet}$, $\maxit$, tol}  
    \State Input: $A$,$\bv$,$\xv_0$ : matrix, right-hand-side, and initial guess; 
    $W\in\Real^{n\times \Na}$ : matrix of $N_{\DeflationSet}$ augmented vectors
    \State $\bv_0 = \bv - A \xv_0$ \Comment Adjusted $\bv$ for initial guess $\xv_0$
    \State $\vv_{\Na+1} = \bv_0$ \Comment Will hold normalized $(I-V_p V_p^T )\bv_0$
    \blue \If{$\Na>0$}
      \State // $V_{1:\Na}$ holds columns $1:\Na$ of $V$.
      \State $A W = V_{1:\Na} \, H_{1:\Na,1:\Na}$ \Comment Compute QR factorization of $A W$.  \label{eq:agmresQR}
        \State $\vv_{\Na+1} = \vv_{\Na+1} - V_p (V_p^T \vv_{\Na+1})$  \Comment $\vv_{\Na+1} = (I-V_p V_p^T )\bv$
    \EndIf \black
    \State $\vv_{\Na+1} = \vv_{\Na+1}/\| \vv_{\Na+1} \|_2$ \Comment $\vv_{\Na+1}$ is column $\Na+1$ of $V$
    \State $Q=I_{1:m,1:\Na+1}$ \Comment Holds product of Givens rotations
    \For{$k=1:\Na+1$}
       \State $g(k) = V(:,k)^T \bv_0$   \Comment $\gv$ is the RHS for the least squares solve
    \EndFor
    \For{$k=\Na+1,\Na+2,\ldots,\Na+\maxit$} \Comment Arnoldi iterations 
      \State $\vv_{k+1} = A \vv_{k}$  \Comment WaveHoltz wave-solve with initial condition $\vv_{k}$.
      \For{$i=1:k$} \Comment Make $\vv_{k+1}$ orthogonal to $\vv_i$
         \State $h_{ik} = \vv_i^T\vv_{k+1}$; \quad $\vv_{k+1} = \vv_{k+1} - h_{ik}\vv_i$
      \EndFor
      \State $h_{k+1,k}= \| \vv_{k+1} \|_2$; \quad $\vv_{k+1}= \vv_{k+1}/\| \vv_{k+1}\|_2$ \Comment Column $k+1$ of $V$
      \State // Apply Givens to H (part of solving final least squares problem):
      \State $H(1:k,k) = Q(1:k,1:k) \, H(1:k,k)$ \Comment Apply previous rotations to new column 
      \State $\rho=H(k,k);$ ~  $H(k,k) = \sqrt{\rho^2 + H(k+1,k)^2};$ 
      \State $c = \rho/H(k,k);$ ~$s = H(k+1,k)/H(k,k);$~ $ H(k+1,k) = 0;$
      \State // Apply Givens rotation to Q 
      \State $Q(k+1,:) = -s\, Q(k,:);$ ~ $Q(k,:) = c\, Q(k,:);$ ~ $Q(k+1,k+1) = c;$ ~ $Q(k,k+1) = s;$
      \State // Apply Givens rotation to RHS g 
      \State $g(k+1) = -s\, (k,1);$ ~ $g(k) = c\, g(k);$  \Comment $g(k+1)$ holds the current 2-norm residual
      \State \textbf{if}~ $|g(k+1)| < \rm{tol} \, \| \bv \|_2$ ~\textbf{then} ~break; ~\textbf{end if}
    \EndFor
    \State // Solve $\yv = \argmin_{\yv \in \Real^{k}} \, \|  V_{1:k+1}^T \bv - H \yv \|_2$ 
    \State $H_{1:k,1:k} \,\yv = \gv_{1:k}$  \Comment Solve the triangular square system \label{eq:agmresFinalSolve}
    \State $\xv = \xv_0  + \begin{bmatrix} \blue  W  \black ~\vertbar~ V_{\Na+1:k}\end{bmatrix} \yv$ \Comment Approximate solution
 \EndFunction
\end{algorithmic} 
\label{alg:augmentedGmresWithGivens}
\end{algorithm}

Given an initial guess, $\xv_0$, the augmented GMRES method will find 
a solution of the form $\xv = \xv_0 + \xvp$
where $\xvp$ is in the
Krylov space generated by $A$ and $\bv_0 \eqdef \bv - A\xv_0$,  augmented with the space spanned by the vectors $\wv_m$, 
\ba
    \xvp \in {\rm span}\{ \wv_1, \wv_2, \ldots \wv_p,\} \cup \, \Kc_j(A,\bv_0) . 
    \label{eq:solutionSpaceAGMRES}
\ea
The solution $\xvp$ is chosen to minimize the two-norm of the residual $\| \bv_0 - A \xvp\|_2$.
The method starts on line~\ref{eq:agmresQR} by forming the reduced QR factorization of $A W$
\ba
  & A W = V_{1:\Na} \, H_{1:\Na,1:\Na}, \qquad V_{1:\Na}\in\Real^{n\times \Na}, \quad H_{1:\Na,1:\Na}\in\Real^{\Na\times \Na} , \\
& V_{1:\Na} =
\begin{bmatrix}
       \vertbar & \vertbar &  & \vertbar \\
       \vv_{1} & \vv_2 & \ldots & \vv_{\Na}    \\
       \vertbar & \vertbar &  & \vertbar \\
       \end{bmatrix},   
  \qquad
    H_{1:\Na,1:\Na} =
    \begin{bmatrix}
      h_{11} & h_{12} & \ldots & h_{1p} \\
             & h_{22} & \ldots & h_{2p} \\
             &        & \ddots & \vdots \\
             &        &        & h_{pp}
    \end{bmatrix} .
\ea
Here the notation $V_{1:\Na}$ indicates that this matrix holds columns $1$ to $\Na$ of matrix $V$.
Note that when the augmented vectors $\wv_m$ are the true eigenvectors with eigenvalues $\lambda_{h,m}$ 
the product $A W$ can be computed efficiently by multiplying column $m$
of $W$ by the corresponding eigenvalue of $A$, $\mu_{h,m} = 1-\beta(\lambda_{h,m})$. Otherwise $A\vv_m$ involves a wave-solve (which could be pre-computed for efficiency).
Let $\bv_0'$ be the projected value of $\bv_0$ that removes the components of $\bv_0$ along the augmented vectors $\wv_i$, 
$
   \bv_0' = ( I - V_{1:\Na} V_{1:\Na}^T) \, \bv_0 .
$
Column $\vv_{\Na+1}$ of $V$ is $\bv_0'$ normalized by its length 
$
   \vv_{\Na+1} = \bv_0'/\| \bv_0' \|_2.
$
Starting from $\vv_{\Na+1}$, the Arnoldi algorithm is used to generate $j$ additional columns $\vv_{\Na+2},\ldots,\vv_{\Na+j+1}$, and
the remaining entries in the upper Hessenberg matrix 
$H\in\Real^{\Na+j+1\times \Na+j}$ 
\ba
H =
    \begin{bmatrix}
      h_{11} & h_{1,2} & \ldots & h_{1,\Na} &  h_{1,\Na+1} & h_{1,\Na+2} & \ldots & h_{1,\Na+j+1} \\
             & h_{2,2} & \ldots & h_{2,\Na} &   h_{2,\Na+1} & h_{2,\Na+2} & \ldots & h_{2,\Na+j+1} \\
             &         & \vdots &         &    \vdots  &  \vdots   &        &  \vdots   \\
             &         &        &h_{\Na \Na}   &    \vdots  &  \vdots   &        &  \vdots   \\
             &         &        &         &  h_{\Na+1,\Na+1} & h_{\Na+1 \Na+2} & \ldots & h_{\Na+1,\Na+j+1} \\
             &         &        &         &  h_{\Na+2,\Na+1} & h_{\Na+2,\Na+2} & \ldots & h_{\Na+2,\Na+j+1} \\
             &         &        &         &              & h_{\Na+2,\Na+3} & \ldots & h_{\Na+3,\Na+j+1} \\
             &         &        &         &              &            & \ddots & \vdots \\
             &         &        &         &              &            &        & h_{\Na+j+1,\Na+j}
    \end{bmatrix} ,
\ea
such that
\ba
   A  \begin{bmatrix}  W  ~\vertbar~ V_{\Na+1:\Na+j} \end{bmatrix}  
   & = \begin{bmatrix}  A W  ~\vertbar~ A V_{\Na+1:\Na+j} \end{bmatrix}
   = \begin{bmatrix}  V_{1:\Na} \, H_{1:\Na,1:\Na} ~\vertbar~ A V_{\Na+1:\Na+j} \end{bmatrix} ,\\
  & = V_{1:\Na+j+1} H  .
\ea
The AGMRES solution will be a linear combination of the columns of $W$ and the columns of $V_{\Na+1:\Na+j}$ (note that we use $W$, not $V_{1:\Na}$)
\ba
    \xvp = \begin{bmatrix}  W  ~\vertbar~ V_{\Na+1:\Na+j} \end{bmatrix} \yv . \label{eq:gmresSolution}
\ea
Multiplying~\eqref{eq:gmresSolution} by $A$ gives 
\ba
   A \xvp & = \begin{bmatrix}  A W                  ~\vertbar~ A V_{\Na+1:\Na+j} \end{bmatrix} \yv 
            = V_{1:\Na+j+1} H \yv .
\ea
The residual can be written as 
$
   \rv = \bv_0 - A \xvp = \bv_0 - V_{1:\Na+j+1} H \yv 
$
where $\yv$ solves the least squares problem
\ba
  \yv = \argmin_{\yv \in \Real^{\Na+j}} \, \|  \bv_0 - V_{1:\Na+j+1} H \yv \| .  \label{eq:GMRESleastSquareI}
\ea
Multiplying the expression inside the norm in~\eqref{eq:GMRESleastSquareI} by $V_{1:\Na+j+1}^T$ gives the final form of the least squares problem
\ba
  \yv = \argmin_{\yv \in \Real^{\Na+j}} \, \|  V_{1:\Na+j+1}^T \bv - H \yv \| .
\ea
This solution for $\yv$ is found on line~\ref{eq:agmresFinalSolve} where, due to the applications of the Givens rotations,
the matrix $H_{1:k,1:k}$ is upper triangular.
After solving for $\yv$, $\xvp$ is found from~\eqref{eq:gmresSolution} and then $\xv=\xv_0+\xvp$.

\subsection{Augmented (recycled) BICGSTAB (ABICGSTAB)}  \label{sec:augmentedBICGSTAB}

The augmented (recycled) BICGSTAB algorithm (ABICGSTAB) from Amritkar et alia~\cite{Amritkar2015} is given in Algorithm~\ref{alg:augmentedBiCGStab}.
The algorithm in~\cite{Amritkar2015} is called the recycled BICGSTAB method where recycling refers to the technique of re-using vectors from 
previous BICGSTAB computations to augment the space.
The standard BICGSTAB method can be used for non-symmetric systems and like standard CG uses a short recurrence with fixed storage requirements.
Standard BICGSTAB will thus generally use fewer floating point operations and less storage than standard GMRES.
The augmented version ABIGSTAB will involve more operations and 
require additional storage associated with the $\Na$ augmented vectors.
Some of the additional costs can be offset by precomputing some quantities that can be reused for different right-hand-sides.

\renewcommand{\algFontSize}{\small}
\begin{algorithm}
\algFontSize 
\caption{Augmented (recycled) BiCGStab, from Amritkar et alia~\cite{Amritkar2015}.}
\begin{algorithmic}[1]

  \Function{$\xv_j=$augmentedBiCGStab}{$A$,$\bv$,$W$,$\xv_0$, ${\rm tol}$, $\maxit$} 
    \State Input: $A$,$\bv$,$\xv_0$ : matrix, right-hand-side, and initial guess; 
       $W\in\Real^{n\times \Na}$ : matrix of $\Na$ augmented vectors
    \State $\xv=\xv_0$, $\rv = \bv - A \xv $
    \State $\rvTilde = \rv$   \Comment Arbitrary vector with $\rvTilde^T\rv \ne 0$
    \State \blue $A W = Q R $  \Comment Thin QR decomposition of $A W$.  
    \State $\etav_1 = Q^T \rvHat$,  $\rv = \rv - Q \etav_1 \rv$ \Comment $\rv \leftarrow (I-Q Q^T) \rv$
    \State $\xiv=-\etav_1$ \Comment $\xiv$ holds $W$ updates for $\xv$ used in the line~\ref{alg:xUpdate}  \black
    \For{$i=1,2,\ldots,\maxit$}   \label{eq:bicgstabStartIterations}
       \State $\rho = \rvTilde^T \rv$ \Comment If $\rho=0$, breakdown occurs, exit gracefully.
       \State \textbf{if} $i\equiv 0$ \textbf{then} $\pv = \rv$ \textbf{else} $\beta=(\rho/\rhoOld)(\alpha/\omega)$; $\pv = \rv + \beta ( \pv-\omega \vv)$ \textbf{end}
       \State $\vv = A \pv$  \Comment Wave-solve.
       \State \blue $\etav_1=Q^T \vv$; $\vv=\vv-Q\etav_1$ \Comment $\vv \leftarrow (I-Q Q^T) \vv$ \black  \label{eq:abicgstabProjectV}
       \State $\alpha =\rho/(\rvTilde^T \rv)$; $\sv=\rv-\alpha \vv$
       \If{ $\| \sv\|_2 \le {\rm tol}\, \| \bv \|_2$ }
         \State $\xv = \xv + \alpha \pv$; $\rv=\sv$; \blue $\xiv=\xiv+\alpha \etav_1$;  \black\textbf{break}
       \EndIf
       \State $\tv = A \sv$   \Comment Wave-solve.
       \State \blue $\etav_2=Q^T \tv$; $\tv=\tv - Q\etav_2$  \Comment $\tv \leftarrow (I-Q Q^T) \tv$ \black \label{eq:abicgstabProjectT}
       \State $\omega= (\tv^T\sv)/(\tv^T\tv)$; \blue $\xiv=\xiv + \alpha \etav_1 + \omega \etav_2$; \black $\xv=\xv + \alpha \pv + \omega \sv$;  $\rv = \sv-\omega \tv$ 
       \State \textbf{If} $\| \rv\|_2 \le {\rm tol} \,\| \bv\|_2$ \textbf{then}  \textbf{break}
       \State $\rhoOld=\rho$
    \EndFor
    \State $\xv = \xv \blue - W (R^{-1} \xiv)$ \label{alg:xUpdate} \black \Comment $R$ is upper triangular.
 \EndFunction
\end{algorithmic} 
\label{alg:augmentedBiCGStab}
\end{algorithm}

\mni
The changes to the standard BICGSTAB algorithm are highlighted in {\blue blue} in Algorithm~\ref{alg:augmentedBiCGStab}.
As with AGMRES, ABICGSTAB finds a solution $\xv$ from the augmented space in~\eqref{eq:solutionSpaceAGMRES}.
Also similarly to AGMRES, ABICGSTAB computes the reduced $QR$ factorization of $AW$ 
and uses $Q$ to project the initial conditions and later iterates.
See the comments in Section~\ref{sec:augmentedGMRES} about efficient ways to compute $AW$ and precomputing the QR factorization.
The main iterations start on line~\ref{eq:bicgstabStartIterations} and there are up to two matrix-vector products per iteration.
The main additional cost per iteration is in the projection steps on lines~\ref{eq:abicgstabProjectV} and~\ref{eq:abicgstabProjectT}.



\section{Numerical experiments on single grids using summation by parts operators\label{sec:numerical}}

For the examples presented in this section we consider the wave equation with constant  speed of sound equal to one and with homogeneous Dirichlet boundary conditions discretized on a single curvilinear grid. The curvilinear grid is described by a mapping 
{
\[
   (x(r,s),y(r,s)) = \Gv(r,s) , \ \ (r,s) \in [-1,1]^2,
\]
where $\Gv$ is assumed to be smooth and invertible.
{The parameter space $(r,s)$ is discretized with a Cartesian grid $r_i=-1+i\dr$, $i=0,1,2,\ldots,N_r$, and $s_j=-1+j\ds$, $j=0,1,2,\ldots,N_s$ where $\dr=2/N_r$ and $\ds=2/N_s$.}
}
The Laplacian is discretized using the conservative formulation described, for example, in \cite{AppPet09}. The metric coefficients and derivatives are approximated using summation by parts operators introduced in \cite{MatNor04,Mattsson2012} and as provided by the Julia package \cite{ranocha2021sbp}. The boundary conditions are enforced by projection (see for example \cite{Eriksson1855469}).

When using SBP operators, the discrete approximation to the Laplacian, $\Lph$, is in general not symmetric but can be symmetrized by a diagonal and positive definite matrix $M_{p,h}$ (usually referred to as the mass matrix). Most of this matrix is a scalar multiple of the identity matrix, only a few diagonal entries at the top left and bottom right,  corresponding to grid points near the boundary, are different than the constant value in the interior (see e.g. \cite{MatNor04}).    

Thus when using the implicit time stepping scheme requires us to compute solutions to non-symmetric systems of equations on the form 
\[
\left(I - \frac{\Delta t^2}{2} \Lph \right) (W^{n+1}_\jv+W^{n-1}_\jv) = g.
\]
Here $g$ contains the solution at time $t^n$ and forcing terms. In order to use PCG, preconditioned with (A)MG it is preferable to first symmetrize by multiplying with $M_{p,h}$ and instead solve the symmetric system of equations
\ba
\left(M_{p,h} - \frac{\Delta t^2}{2} M_{p,h} \Lph \right) (W^{n+1}_\jv+W^{n-1}_\jv) = M_{p,h} g. \label{eq:symm_imp1}
\ea
Note that the structure of this system of linear equations is identical to the system of linear equations that would result from a symmetric interior penalty discontinuous Galerkin or a continuous Galerkin formulation. 

Now, the {WaveHoltz matrix $A = I - S_h$} is a function of $\Lph$ and is therefore no longer symmetric. It can also be symmetrized by left multiplication with $M_{p,h}$. Doing so we find that one option for solving Helmholtz equations by the WaveHoltz method then amounts to solving the symmetric positive definite linear system of equations       
\ba
    M_{p,h}(I - S_h) \Uv = M_{p,h}\bv.  
\ea

Alternatively, and this is the approach we take here, we can introduce $Z^{n} = M_{p,h}^{\frac{1}{2}} W^{n}$ and evolve this variable using the same scheme   
\ba
\left(I - \frac{\Delta t^2}{2} M_{p,h}^{\frac{1}{2}} \Lph M_{p,h}^{-\frac{1}{2}} \right) (Z^{n+1}_\jv+Z^{n-1}_\jv) = 2Z^{n}_\jv -  \Delta t^2 M_{p,h}^{\frac{1}{2}} f(\xv_\jv) \, \cos(\omega t^n). \label{eq:symm_imp2}
\ea
Since the matrix $M_{p,h}^{\frac{1}{2}} \Lph M_{p,h}^{-\frac{1}{2}}$ shares eigenvalues with $\Lph$ we can directly apply the deflation of the forcing to the symmetric problem       
\ba \label{eq:sqrt_lin_sys}
    (I - S_h) \Zv = M_{p,h}^{\frac{1}{2}} \bv,  
\ea
where the action of $S_h$ and $\bv = \Pi_h 0$ are computed using (\ref{eq:symm_imp2}). The deflation of the forcing is then done to $M_{p,h}^{\frac{1}{2}} f(\xv_\jv)$ using the orthogonal eigenvectors of $M_{p,h}^{\frac{1}{2}} \Lph M_{p,h}^{\frac{1}{2}}$. Once the system (\ref{eq:sqrt_lin_sys}) has been solved, the solution to Helmholtz equation is obtained as $\Uv = M_{p,h}^{-\frac{1}{2}} \Zv$. Note that since $M_{p,h}$ is diagonal, operations with $M_{p,h}^{-\frac{1}{2}}$ or $M_{p,h}^{\frac{1}{2}}$ are trivial. 

\medskip
{The EigenWave algorithm~\cite{EigenWave} is used for the computations of the eigenpairs used in this section.
EigenWave is built on top of a matrix-free Arnoldi method and here we use 
{\tt ArnoldiMethod.jl}, a Julia implementation of the Arnoldi method with a Krylov-Schur restart.
Each matrix-vector product in the Arnoldi algorithm corresponds to application of the WaveHoltz matrix $A$, which in turn involves the solution to the wave equation.
}


\subsection{Comparison with standard CG} \label{seq:SBPex1}

In a first example we consider the domain 
\[
 x(r,s) = r + 0.1\sin(2s), \ \  y(r,s) = s - 0.3 \sin(2r).
\]
We use fourth order accurate SBP operators and determine the number of grid points according 
to the rule of thumb~\eqref{eq:NxFromPPW} with $\epsilon = 10^{-4}$, and $\domainSize=2.25$. We use WaveHoltz with $N_p=3$ periods and $10$ timesteps per period. 
The computations presented in this section all use the Gaussian source term (\ref{eq:gaussianSource}) with $ x_0 = 0.1, y_0= 0.2$ and $\alphag=\omega^2$, and $\betag=\omega$.

{
\newcommand{\figw}{0.40\textwidth}
\newcommand{\figh}{0.33\textwidth}
\begin{figure}[htb]
\begin{center}
\begin{tikzpicture}
  \useasboundingbox (0,.8) rectangle (2.1*\figw,\figh);  
  \begin{scope}[yshift=0*\figh]
    \figByWidth{0.00*\figw}{0*\figh}{fig/solution_om20_nev50}{\figw}[0.11][0.][0.][0.]
    \figByWidth{1.05*\figw}{0*\figh}{fig/solution_om40_nev100}{\figw}[0.11][0.][0.][0.]
  \end{scope}  
\end{tikzpicture}
\end{center}
\caption{Solutions to the two problems considered in Section~\ref{seq:SBPex1}. 
   Left: $\omega=20$. Right: $\omega=40$.
\label{fig:SBP_example_1_sol}}
\end{figure}
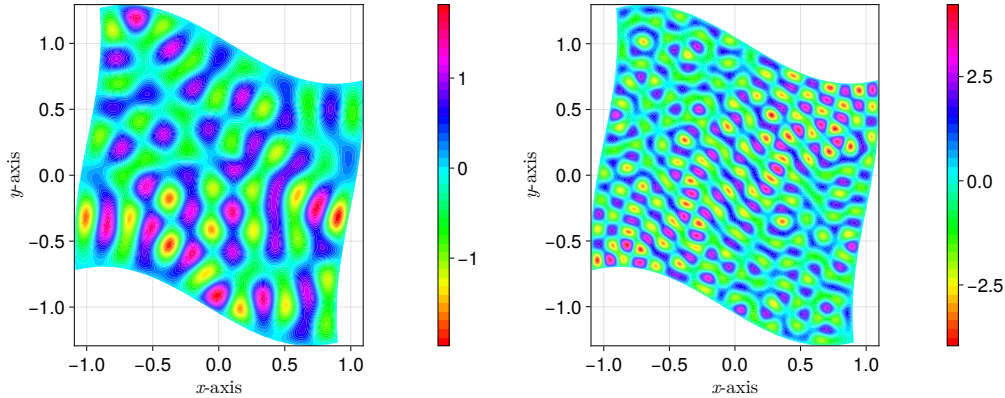
}
 
We first consider $\omega = 20$ which results in a grid with $318 \times 318$ unknowns (we display the solution to this problem in Figure \ref{fig:SBP_example_1_sol}).
We use $\Ndeflate=50$ eigenpairs for the deflation. {These are computed using EigenWave} to a tolerance of $10^{-14}$ (the number of wave solves to reach this tolerance is $\Nwseig=150$ wave-solves or 3 solves per eigen pair). The location of the deflated eigenvalues on the $\beta$-curve are displayed  in Figure \ref{fig:SBP_example_1_bet} where the ACR can also be read off. 

{
\newcommand{\figw}{0.40\textwidth}
\newcommand{\figh}{0.35\textwidth}
\begin{figure}[htb]
\begin{center}
\begin{tikzpicture}
  \useasboundingbox (0,.75) rectangle (2.1*\figw,\figh);  
  \begin{scope}[xshift=0*\figh]
    \figByWidth{0}{0}{fig/residual_om20_nev50}{\figw}[0.11][0.][0.][0.]
    \draw[black] (3.25,5.65) node {$\omega=20$, $\Ndeflate=50$ eigenvectors deflated};
  \end{scope}
  \begin{scope}[xshift=1.05*\figw]
    \figByWidth{0}{0}{fig/residual_om40_nev100}{\figw}[0.11][0.][0.][0.]
    \draw[black] (3.25,5.65) node {$\omega=40$, $\Ndeflate=100$ eigenvectors deflated};
  \end{scope}  
\end{tikzpicture}
\end{center}
\caption{Residuals for the various methods used to solve the two problems considered in Section \ref{seq:SBPex1}. 
  Left: $\omega=20$ with $\Ndeflate=50$ eigenvectors deflated.
  Right: $\omega=40$ with $\Ndeflate=100$ eigenvectors deflated.
  DCG: WaveHoltz and the deflated conjugate gradient algorithm.
  CG: WaveHoltz accelerated with the conjugate gradient method (no deflation).
  CG with $f_d$: WaveHoltz and conjugate gradients with forcing deflated.
  WHI with $f_d$: WaveHoltz fixed-point iteration with forcing deflated.
  ACR line: theoretical asymptotic convergence rate for the FPI with deflation (see Figure~\ref{fig:SBP_example_1_bet}).
\label{fig:SBP_example_1_res}}
\end{figure}
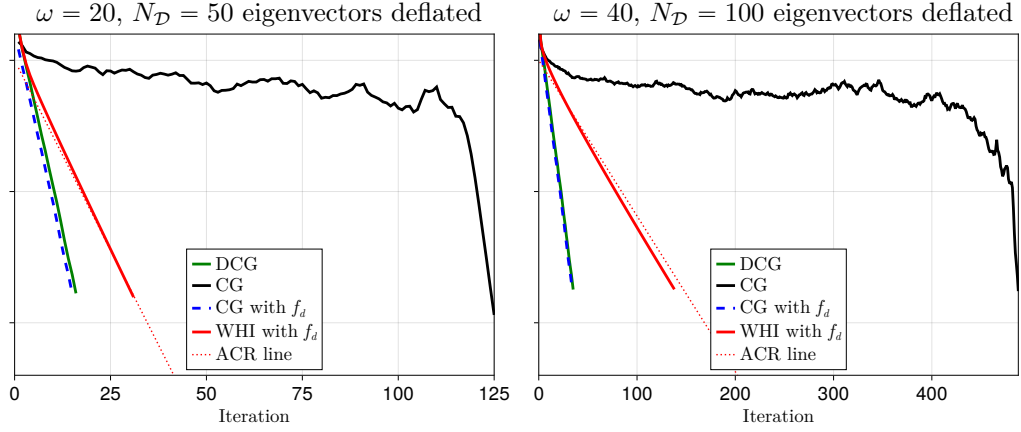
}

Figure  \ref{fig:SBP_example_1_res} displays  the residuals for this problem when using the standard CG accelerated WaveHoltz with with and without the forcing being deflated (denoted ``CG'' and ``CG with $f_d$''), the plain WaveHoltz iteration applied to the deflated forcing (denoted ``WHI with $f_d$''), and the DCG method (denoted ``DCG'') applied to the problem with the non-deflated forcing. As can be seen CG with $f_d$ and DCG have the lowest iteration count, followed by WHI with $f_d$ (which rate of convergence is well predicted by the ACR). The standard CG method is significantly slower.        

We then repeat with $\omega = 40$ on a grid with $756 \times 756$ unknowns and with $\Ndeflate=100$ eigenpairs for the deflation. 
This time the $100$ eigenpairs pairs are computed with {the EigenWave algorithm} to a tolerance of $10^{-14}$ using at total of $\Nwseig=593$ wave-solves. This is slightly less efficient and can be explained from the location of the deflated eigenvalues on the $\beta$-curve are displayed to the right in Figure \ref{fig:SBP_example_1_bet} (see also \cite{EigenWave}). 

{
\newcommand{\figw}{0.40\textwidth}
\newcommand{\figh}{0.30\textwidth}
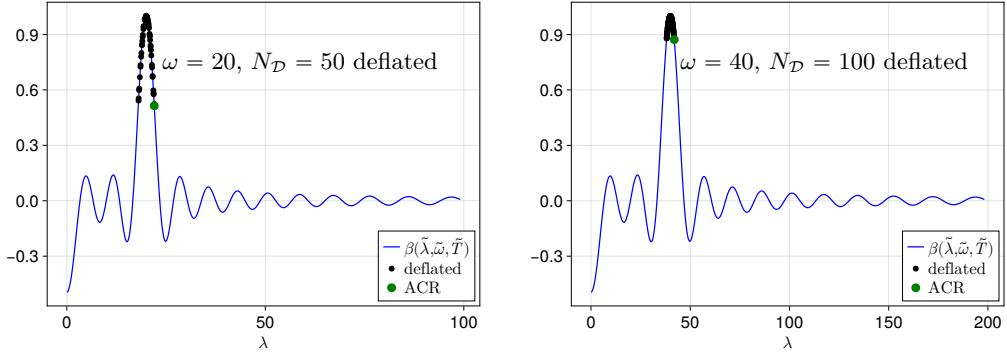
\begin{figure}[htb]
\begin{center}
\begin{tikzpicture}
  \useasboundingbox (0,.75) rectangle (2.1*\figw,\figh);  
  \begin{scope}[xshift=0*\figh]
    \figByWidth{0}{0}{fig/beta_om20_nev50}{\figw}[0.][0.][0.][0.]
    \draw[black] (4,4) node {\small $\omega=20$, $\Ndeflate=50$ deflated};
  \end{scope}
  \begin{scope}[xshift=1.05*\figw]
    \figByWidth{0}{0}{fig/beta_om40_nev100}{\figw}[0.][0.][0.][0.]
    \draw[black] (4,4) node {\small $\omega=40$, $\Ndeflate=100$ deflated};
  \end{scope}  
\end{tikzpicture}
\end{center}
\caption{Beta functions and deflated eigenvalues for the problems described in Section~\ref{seq:SBPex1} and shown in Figure~\ref{fig:SBP_example_1_res}.
\label{fig:SBP_example_1_bet}}
\end{figure}
}

The right of Figure \ref{fig:SBP_example_1_res} again displays the residual. As can be seen CG with $f_d$ and DCG again have the lowest iteration count, followed by WHI with $f_d$ (which rate of convergence is well predicted by the ACR). Again, the standard CG method is significantly slower.     

Ignoring the cost of the deflation and inflation we can estimate the cost to solve for $N_f$ right hand sides. If, for the lower frequency we take the number of wave solves to reach convergence to be $\approx 20$ for the deflated (D)CG method and $\approx 125$ for the standard CG method we can estimate the cost (in number of wave solves) to be 
\bse
\ba
\text{Cost deflated WaveHoltz} \approx 20 \, N_f + \Nwseig &\approx 20 \, N_f + 150,  \\  
 \text{Cost standard WaveHoltz} &\approx 125 N_f. 
\ea 
\ese
From this we can conclude that (in terms of number of wave solves) the deflated method can be faster already with as few as two right hand sides. The same computation for the $\omega = 40$ case yields the same conclusion.

\subsection{Varying the number of eigenvectors with frequency} \label{sec:SBPex2}

From Figure \ref{fig:SBP_example_1_bet} we see that, even if we increase the number of eigenvectors in the deflation linearly with increasing $\omega$, the ACR becomes worse. 
{In order to obtain frequency independent convergence rates for increasing frequency it is necessary to increase the number of deflated eigenvectors, $\Ndeflate$,
in proportion to $\omega^d$, in $d$ dimensions, as will now be shown. 
From the form of the beta function in~\eqref{eq:betaGap}, it can be seen that we need to deflate the eigenvectors with eigenvalues $\lambda_{h,m}$ that satisfy
\ba
  \left| \f{\lambda_{h,m}-\omega}{\omega} \right| <  C_\lambda, \label{eq:freqIndepence}
\ea
where $C_\lambda$ is some positive constant. 
Using~\eqref{eq:freqIndepence} in~\eqref{eq:betaGap} implies that with deflation the asymptotic convergence rate will be approximately
\ba
   \ACR \approx 1 - C_\beta C_\lambda^2.
\ea
According to Weyl's asymptotic analysis for the distribution of eigenvalues for the Laplace operator with Dirichlet boundary conditions, the number
of eigenvalues per unit length near $\lambda\approx \omega$ scales as $\omega^{d-1}$ in $d$ space dimensions. For example in one dimension the eigenvalues are equally spaced for all $\omega$.
Combining~\eqref{eq:freqIndepence} with Weyl's result, implies that \textbf{to achieve frequency independent convergence rates the number of deflated eigenvalues should scale as $\Ndeflate = O(\omega^d)$ in $d$ space dimensions}.
}

\medskip
To study this result numerically in some detail we compute solutions to the problem described in Section \ref{seq:SBPex1} for $\omega = 5,6,\ldots,29,30,$ and with all parameters 
equal except for taking  $\epsilon = 10^{-3}$ {in the PPW rule of thumb}. We chose the number of eigenvectors to grow either linearly with $\omega$ as $\Ndeflate = \big\lfloor \frac{3\omega}{2} \big\rfloor$, 
or growing quadratically with $\omega$ as $\Ndeflate = \big\lfloor \frac{\omega^2}{10} \big\rfloor$. 

{
\newcommand{\figw}{0.40\textwidth}
\newcommand{\figh}{0.30\textwidth}
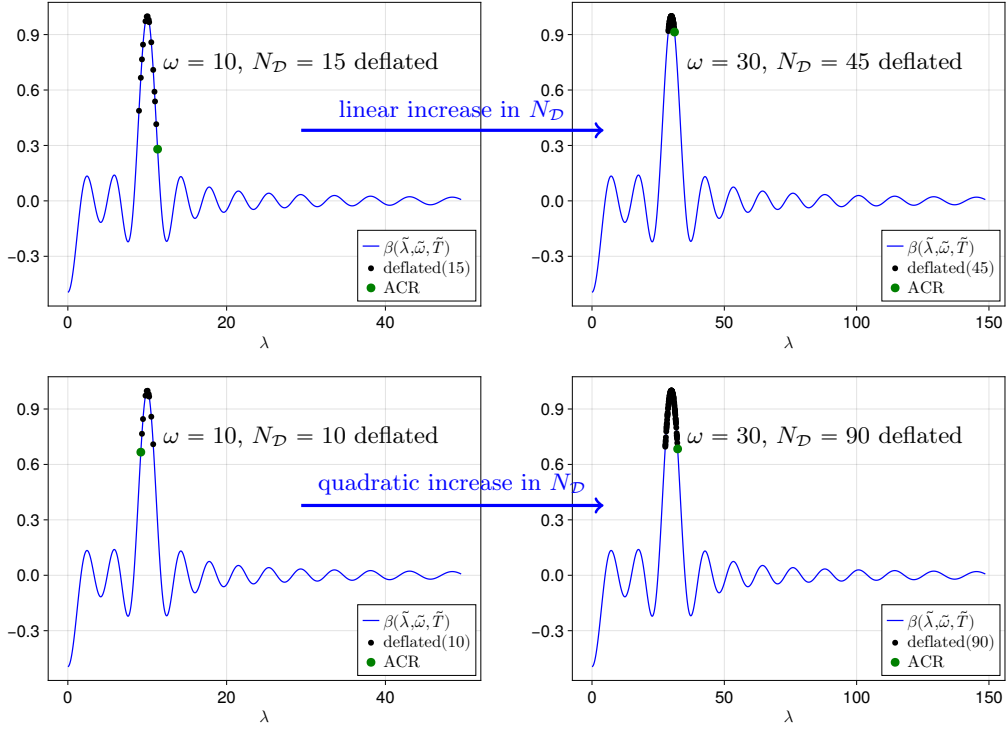
\begin{figure}[htb]
\begin{center}
\begin{tikzpicture}
  \useasboundingbox (0,.8) rectangle (2.1*\figw,2*\figh);  
  \begin{scope}[yshift=1*\figh]
    \begin{scope}[xshift=0*\figw]
      \figByWidth{0}{0}{fig/linear_beta_om_010_nev015}{\figw}[0.][0.][0.][0.]
      \draw[black] (4,4) node {\small $\omega=10$, $\Ndeflate=15$ deflated};
    \end{scope}
    \begin{scope}[xshift=1.05*\figw]
      \figByWidth{0}{0}{fig/linear_beta_om_030_nev045}{\figw}[0.][0.][0.][0.]
      \draw[black] (4,4) node {\small $\omega=30$, $\Ndeflate=45$ deflated};
    \end{scope}  
    \draw[->,black,very thick,blue] (4,3.11) -- (8,3.11) node[midway,above] {\small linear increase in $\Ndeflate$};
  \end{scope} 

  \begin{scope}[xshift=0*\figh]
    \figByWidth{0}{0}{fig/quadratic_beta_om_010_nev010}{\figw}[0.][0.][0.][0.]
    \draw[black] (4,4) node {\small $\omega=10$, $\Ndeflate=10$ deflated};
  \end{scope}
  \begin{scope}[xshift=1.05*\figw]
    \figByWidth{0}{0}{fig/quadratic_beta_om_030_nev090}{\figw}[0.][0.][0.][0.]
    \draw[black] (4,4) node {\small $\omega=30$, $\Ndeflate=90$ deflated};
  \end{scope}  
  \draw[->,black,very thick,blue] (4,3.10) -- (8,3.10) node[midway,above] {\small quadratic increase in $\Ndeflate$};

\end{tikzpicture}
\end{center}
\caption{
Locations of the deflated eigenvalues on the beta functions.
Top: When the number of deflated eigenvectors increases linearly from left to right, the ACR is seen to increase.
Bottom: When the number of deflated eigenvectors increases quadratically from left to right, the ACR is seen to be nearly constant.
\label{fig:SBP_example_2a}}
\end{figure}
}

In Figure \ref{fig:SBP_example_2a} we display where the computed eigenvalues are located on the $\beta$-function for $\omega = 10 $ and 30 and for the linear and quadratic choices for the number of eigenvalues. As can be seen the ACR is roughly constant for the latter (quadratic) choice.  

{
\newcommand{\figw}{0.32\textwidth}
\newcommand{\figh}{0.25\textwidth}
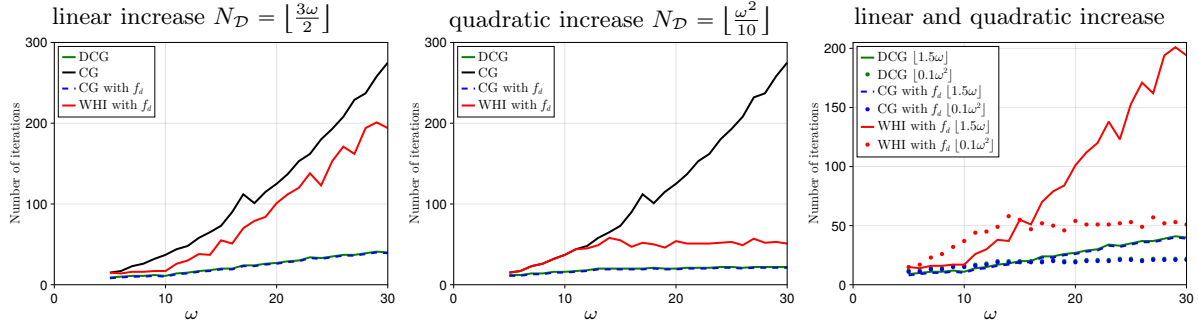
\begin{figure}[htb]
\begin{center}
\begin{tikzpicture}
  \useasboundingbox (0,.45) rectangle (3.1*\figw,\figh);  
  \begin{scope}[xshift=0*\figh]
    \figByWidth{0}{0}{fig/iter_vs_omega_linear}{\figw}[0.][0.][0.08][0.]
    \draw[black] (2.5,3.85) node {\small linear increase $\Ndeflate=\big\lfloor \frac{3\omega}{2} \big\rfloor$};
    \draw[black] (2.5,0) node[yshift=-2pt] {\scriptsize $\omega$};
  \end{scope}
  \begin{scope}[xshift=1.0*\figw]
    \figByWidth{0}{0}{fig/iter_vs_omega_quadratic}{\figw}[0.][0.][0.08][0.]
    \draw[black] (2.75,3.85) node {\small quadratic increase $\Ndeflate=\big\lfloor \frac{\omega^2}{10} \big\rfloor$};
    \draw[black] (2.5,0) node[yshift=-2pt] {\scriptsize $\omega$};
  \end{scope}
  \begin{scope}[xshift=2.0*\figw]
    \figByWidth{0}{0}{fig/iter_vs_omega_linear_quadratic}{\figw}[0.][0.][0.08][0.]
    \draw[black] (2.75,3.85) node {\small linear and quadratic increase};
    \draw[black] (2.5,0) node[yshift=-2pt] {\scriptsize $\omega$};
  \end{scope}     
\end{tikzpicture}
\end{center}
\caption{
 Number of iterations as a function frequency $\omega$ compared to WaveHoltz with conjugate gradients and no deflation (CG).
 Left: the number of eigenvectors is  
$\Ndeflate=\big\lfloor \frac{3\omega}{2} \big\rfloor$.
Middle: $\Ndeflate= \big\lfloor \frac{\omega^2}{10} \big\rfloor$. 
Right: Standard CG results have been removed for easier comparison.
\label{fig:SBP_example_2}}
\end{figure}
}

We also display the number of iterations as a function of $\omega$ for the two choices in Figure \ref{fig:SBP_example_2}.  The number of iterations for the quadratic scaling appears to be approach a constant for all the deflated methods, with the constant being smaller for the CG accelerated methods.

\subsection{Experiment with compressed eigenvectors}

We have seen that increasing the number of eigenvectors in the deflation space reduces the number of iterations. The obvious drawbacks are of course that these vectors must be stored and that it becomes more expensive to carry out the deflation and inflation steps in Algorithm \ref{alg:waveHoltzDirectDeflation}. Compression of the eigenvectors can alleviate the storage requirement but there is a risk that compression, if used too aggressively, can affect the convergence of the iterative method. When working 
on logically Cartesian meshes the singular value decomposition (SVD) is a natural choice for the compression of the eigenvectors. The efficiency of the compression is then measured by how fast the singular values of the eigenvectors decay. In this experiment we numerically investigate this for the setup described for the problem in Section \ref{seq:SBPex1}. 
{Note that throughout this section we {\it do not} need to deflate the iterates, 
the optional deflation step given in line~\ref{eq:deflationIterate} of Algorithm~\ref{alg:waveHoltzDirectDeflation}}.  

Let {$Q \in \Real^{m\times n}$} be the two dimensional matrix representation of an eigenvector {$\Phi_{i,j}$, $i=1,2,\ldots,m$, $j=1,2,\ldots,n$,} 
on a $m \times n$ curvilinear mesh and let 
\[
   Q = USV^T = \sum_{k = 1}^{\min(m,n)} {\bf u}_k \sigma_k {\bf v}_k^T,
\]
be its singular value decomposition. Given a tolerance $\tau_\sigma$, we truncate the SVD so that $r$ is the largest value for which the inequality 
\begin{equation} \label{eq:SVD_trunc}
\| Q - \sum_{k = 1}^{r} {\bf u}_k \sigma_k {\bf v}_k^T \|^2_{F} = \| \sum_{k = r+1}^{\min(m,n)} {\bf u}_k \sigma_k {\bf v}_k^T \|^2_{F} < \tau_\sigma^2,
\end{equation}     
holds. Clearly, if the singular values decay rapidly the compression ratio 
\ba
  \text{compression ratio} = \f{\text{storage for truncated SVD}}{\text{storage for $\Phi_{i,j}$}} = \frac{r(n+m+1)}{nm}
\ea
will be small, signaling a good compression. 

{
\newcommand{\figw}{0.375\textwidth}
\newcommand{\figh}{0.29\textwidth}
\begin{figure}[htb]
\begin{center}
\begin{tikzpicture}
  \useasboundingbox (0.1,.75) rectangle (\figh+2*\figw,\figh);  
  \begin{scope}[xshift=0*\figh]
    \figByHeight{0}{0}{fig/eigenmode1}{\figh}[0.15][0.15][0.0][0.]
    \draw[black] (2.4,4.8) node {\footnotesize Eigenvector, $\lambda_{h,m} \approx \omega$};
  \end{scope}
  \begin{scope}[xshift=.9*\figh]
    \figByWidth{0}{0}{fig/sing_vals_eps}{\figw}[0.][0.][0.0][0.]
    \draw[black] (3.3,4.8) node {\footnotesize Singular values as grid is refined};
  \end{scope}
  \begin{scope}[xshift=.9*\figh+1.0*\figw]
    \figByWidth{0}{0}{fig/ranks_for_eps}{\figw}[0.][0.][0.0][0.]
    \draw[black] (3.0,4.8) node {\footnotesize Estimated rank};
  \end{scope}     
\end{tikzpicture}
\end{center}
\caption{Comparing the SVD compressed and full rank representation of the deflation space. 
   Left: eigenvector $\Phi_{i,j}$ being compressed. 
   Middle: Singular values for $N_x=101$, $179$ and $318$ grid points, corresponding to relative error tolerances of $\eps=10^{-2}$, $10^{-3}$ and $10^{-4}$, respectively.
   Right: estimated rank using an SVD cutoff of $\tau_\sigma = 10^{-5}$.
\label{fig:SBP_example_3}}
\end{figure}
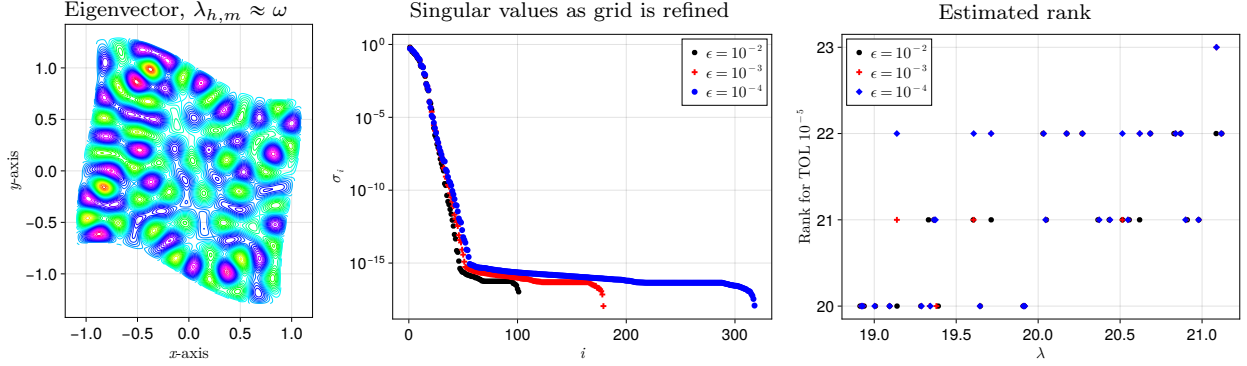
}

To test the ability to compress the eigenfunctions we set $\omega = 20$ and deflate using $\Ndeflate=30$ eigenvectors.
{Grids of different resolutions will be studied, with the number of grid points determined using the PPW rule of thumb~\eqref{eq:NxFromPPW} with relative error tolerances of  
$\epsilon = 10^{-2}, 10^{-3},$ and $10^{-4}$}. 
To the left in Figure \ref{fig:SBP_example_3} we display the the eigenvector corresponding to the eigenvalue closest to $\omega = 20$. To the right we display the singular values ordered in decreasing order for this eigenvector for  $\epsilon = 10^{-2}, 10^{-3},$ and $10^{-4}$ (corresponding to $101\times 101$, $179 \times 179$ and $318 \times 318$ grid points). From this figure it is clear that the decay of the singular values is largely independent to the grid spacing and the compression ratio be increasingly favorable as the mesh is refined.

{
\newcommand{\figw}{0.40\textwidth}
\newcommand{\figh}{0.30\textwidth}
\begin{figure}[htb]
\begin{center}
\begin{tikzpicture}
  \useasboundingbox (0,.80) rectangle (2.1*\figw,1.1*\figh);  
  \begin{scope}[xshift=0*\figh]
    \figByWidth{0}{0}{fig/compare_SVD_FULL_iter}{\figw}[0.][0.][0.][0.]
    \draw[black] (3.25,5.1) node {\small Deflated WaveHoltz, iterations versus $\omega$};
  \end{scope}
  \begin{scope}[xshift=1.05*\figw]
    \figByWidth{0}{0}{fig/rank_vs_omega}{\figw}[0.][0.][0.][0.]
    \draw[black] (3.5,5.1) node {\small Estimated rank versus $\omega$};
  \end{scope}  
\end{tikzpicture}
\end{center}
\caption{
  Left: comparing the number of WaveHoltz iterations as a function of $\omega$ when using the SVD compressed and full rank representation of the deflation space. 
  The black lines show results when choosing the number of deflated eigenvectors to be $\Ndeflate=\big\lfloor \frac{3\omega}{2} \big\rfloor$. The red lines show results for $\Ndeflate=\big\lfloor \frac{\omega^2}{10} \big\rfloor$.
   It is seen that the use of the compressed eigenvectors has no effect on the number of iterations.
   Right: Estimated rank as a function of $\omega$, which grows approximately linearly.
\label{fig:SBP_example_3b}
}
\end{figure}
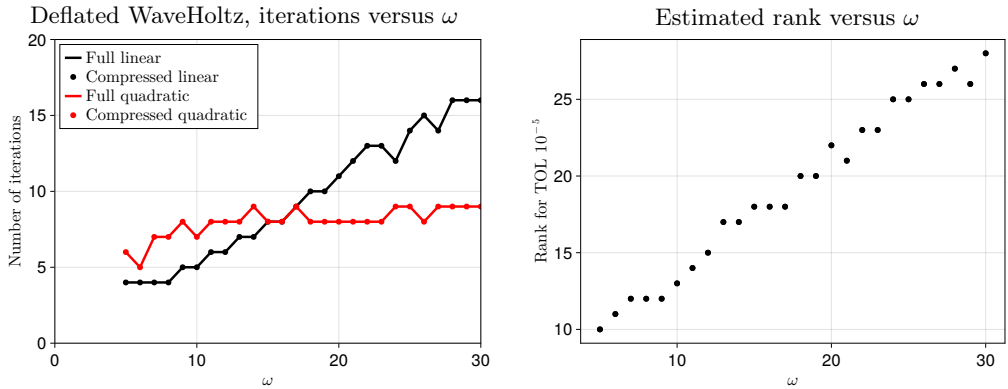
}

We then fix $\tau_\sigma = 10^{-5}$ (smaller than all $\epsilon$) and, for all of the 30 eigenvectors, compute the number $r$ so that the inequality in (\ref{eq:SVD_trunc}) is satisfied. The result is displayed to the right in Figure \ref{fig:SBP_example_3} where we plot $r$ as a function of the eigenvalue corresponding to the eigenvector. As can be seen the rank $r$ is between 20 and 23. Note that the compression ratios for the three deflation spaces are $0.42, 0.24, 0.13$ for $\epsilon = 10^{-2}, 10^{-3},$ and $10^{-4}$, respectively. 

Finally we repeat the experiment from Section \ref{sec:SBPex2} and compute solutions to the problem described in Section \ref{seq:SBPex1} for $\omega = 5,6,\ldots,29,30,$ and with $\epsilon = 10^{-3}$ and $\tau_\sigma = 10^{-5}$. We choose the stopping tolerance in CG (applied to the forcing deflated using the compressed eigenvectors) to be $10^{-4}$. Again we chose the number of eigenvectors to grow either linearly with $\omega$ as $\Ndeflate=\big\lfloor \frac{3\omega}{2} \big\rfloor$, or growing quadratically with $\omega$ as $\Ndeflate=\big\lfloor \frac{\omega^2}{10} \big\rfloor$. In the left part of Figure \ref{fig:SBP_example_3b} we report iteration counts when using the original eigenvectors and when using the compressed eigenvectors. As can be seen the results are identical for the compressed and full rank eigenvectors.  To the right in the same figure we display the rank of the eigenvector whose eigenvalue is closest to resonance as a function of $\omega$. A linear increase with $\omega$ is observed.     

A careful numerical (and perhaps theoretical) study of the interplay between the tolerance used in the compression and the tolerance used in CG is left for the future but we note that we do observe that CG fails to converge (the residual stalls at a certain level) if $\tau_\sigma$ is larger than the {CG tolerance}. {However, even when CG fails to converge the approximate compressed eigenvectors could still potentially be used with the augmented GMRES or BICGSTAB algorithms.}

\subsection{Compressibility of eigenvectors for different meshes}

To examine how the quality of the mesh impacts the compressibility we compute the most slowly converging eigenvector for $\omega = 10$, 
{using different grid resolutions from the PPW rule of thumb~\eqref{eq:NxFromPPW} with relative error tolerances
$\epsilon = 10^{-3}$ and $\epsilon = 10^{-4}$}, on the parametrized grid    
\ba
 x(r,s) = r + a\cos(\frac{\pi s}{2}), \ \  y(r,s) = s + a\cos(\frac{\pi r}{2}).  \label{eq:skewMapping}
\ea
{The parameter $a$ in~\eqref{eq:skewMapping} adjusts the skewness of the grid, with $a=0$ being no skew and $a=0.5$ being highly skewed.}
The singular values (ordered in decreasing order) are plotted for various values of the parameter $a$ in Figure \ref{fig:SBP_example_3c}. In the same figure we have also plotted some of the grids used (down-sampled for visualization purposes). In Figure \ref{fig:SBP_example_3c} we observe that for a Cartesian grid ($a=0$) the rank of the eigenvector is approaching the multiplicity of the eigenvalue. For curvilinear grids we observe a mild growth of the ``rank'' of the eigenvector. Considering that the eigenvectors converge to those of an elliptic problem with variable coefficients it is perhaps expected that such eigenvectors will have low rank as long as the variable coefficients are smooth, as they are on curvilinear and overlapping grids. 

The eigenvector for the curvilinear grid shown in the left plot of Figure~\ref{fig:SBP_example_3} resembles an eigenvector for the square that has been deformed to the curvilinear domain.
One might wonder whether it is possible to avoid the computation of the SVD altogether by expanding the eigenvectors for a curvilinear grid with $a \ne 0$ 
in terms of the known eigenvectors of a Cartesian grid that have been mapped to the curvilinear grid.
On an $m\times m$ grid for the square $[-1,1]^2$ (corresponding to $a=0$), 
the eigenvectors of the discretized Laplacian with Dirichlet boundary conditions are explicitly 
known as a tensor product,  $\Phi = \zv_\mu \, \zv_\nu^T$,  of two one-dimensional eigenvectors $\zv_\mu$ and $\zv_\nu$.
The one dimensional eigenvectors $\zv_k$, $k=1,2,\ldots,m$, have components $z_{k,i}= d_k \sin(\pi k (r_i+1)/2)$, $i=1,\ldots,m$, 
and where $r_i$ are the grid points in the parameter space and $d_k$ is a normalization factor so that $\zv_k^T\zv_l=\delta_{k,l}$.
Let $Z$ be the $m\times m$ matrix with columns $\zv_\mu$, $\mu=1,2,\ldots,m$. 
If $Q$ is the matrix representation for an eigenvector on a curvilinear grid then $Q$ can be written as 
\ba
Q = Z \, C \, Z^T  = \sum_{k=1}^m\sum_{l=1}^m c_{k,l} \, \zv_k \, \zv_l^T, 
\ea
for some matrix $C$. It is straightforward to compute $C$ since the matrix $Z$ is an orthogonal matrix.
For $a=0$, all entries of $C$ are zero except for one, say $c_{\mu,\nu}=1$. For $a \ne 0$, $C$ could have all non-zero entries.
Unfortunately, for the case considered here, $C$ does not appear to be sparse or even to have elements decaying rapidly. In the top right graph in Figure \ref{fig:SBP_example_3c} we have displayed the absolute values of the elements in $C$ (sorted according to magnitude) for grids with $a = 0.0,0.1,\ldots,0.5$ and $\epsilon = 10^{-3}$. 
As can be seen in the figure (note the scale of the horizontal axis), except for the case $a=0$, the decay is very slow.

{
\newcommand{\figwa}{0.30\textwidth}
\newcommand{\figw}{0.30\textwidth}
\newcommand{\figh}{0.28\textwidth}
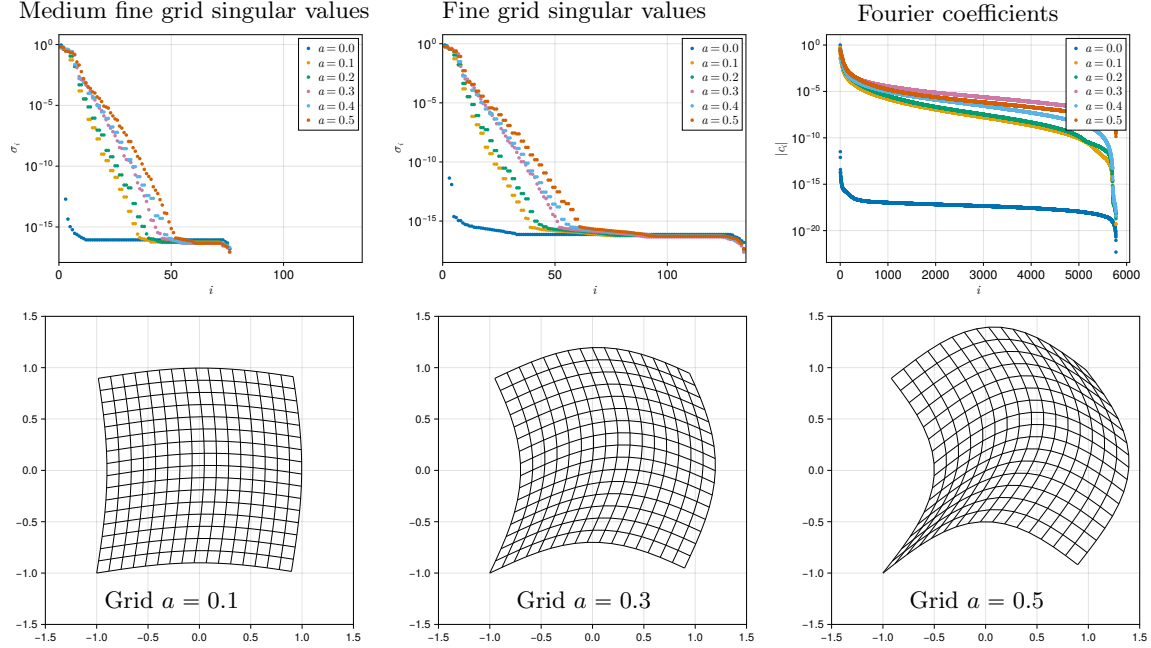
\begin{figure}[htb]
\begin{center}
\begin{tikzpicture}
  \useasboundingbox (0,.60) rectangle (3.1*\figw,1.9*\figh);  

  \begin{scope}[xshift=0*\figh,yshift=\figh]
    \figByWidth{0}{0}{fig/SVDS_of_grid1}{\figwa}[0.][0.][0.][0.]
    \draw[black] (2.55,3.9) node {\small Medium fine grid singular values};
  \end{scope}
  \begin{scope}[xshift=1.025*\figw,yshift=\figh]
    \figByWidth{0}{0}{fig/SVDS_of_grid2}{\figwa}[0.][0.][0.][0.]
    \draw[black] (2.5,3.9) node {\small Fine grid singular values};
  \end{scope} 

  \begin{scope}[xshift=2.05*\figw,yshift=\figh]
    \figByWidth{0}{0}{fig/tensor_proj_coeffs}{\figwa}[0.][0.][0.][0.]
    \draw[black] (2.5,3.9) node {\small Fourier coefficients};
  \end{scope}

 \begin{scope}[xshift=0*\figw]
    \figByWidth{0}{0}{fig/grid_a01}{\figw}[0.1][0.1][0.][0.]
    \draw[black] (2.25,0.75) node {\small Grid $a=0.1$};
  \end{scope}
  \begin{scope}[xshift=1.05*\figw]
    \figByWidth{0}{0}{fig/grid_a03}{\figw}[0.1][0.1][0.][0.]
    \draw[black] (2.5,0.75) node {\small Grid $a=0.3$};
  \end{scope} 
  \begin{scope}[xshift=2.1*\figw]
    \figByWidth{0}{0}{fig/grid_a05}{\figw}[0.1][0.1][0.][0.]
    \draw[black] (2.5,0.75) node {\small Grid $a=0.5$};
  \end{scope} 

  %
\end{tikzpicture}
\end{center}
\caption{Comparing the compressibility of the eigenvector, with eigenvalue close to $\omega=10$,
 on meshes of different quality, measured by the skewness parameter $a$ in the mapping~\eqref{eq:skewMapping}.
  Top left: singular values for different values of $a$ for a medium fine grid.
  Top middle: singular values for different values of $a$ for a fine grid.
  Top right: Coefficients in the expansion of the eigenvector for $a\ne0$ in terms of the eigenvectors for a Cartesian grid $a=0$ that have been mapped to the curvilinear grid.
  Bottom left to right: grids for $a=0.1$, $0.3$ and $0.5$.
\label{fig:SBP_example_3c}}
\end{figure}
}


\section{Numerical experiments on overset grids} \label{sec:oversetGridResults}

In this section we present results from using the WaveHoltz method to solve Helmholtz problems on domains discretized with overset grids.
An overset grid is a collection of overlapping structured grids that covers a domain $\Omega$.
A typical grid consists of one or more Cartesian background grids together with multiple boundary fitted
curvilinear grids. Solution values are matched by interpolation~\cite{CGNS}.
Overset grids allow the use of efficient high-order accurate finite difference schemes for complex geometry with accurate treatment of curved boundaries.
For further details on the implementation of WaveHoltz on overset grids see~\cite{overHoltzPartOne,overHoltzPartTwo}.
The discrete eigenvalues and eigenvectors used in this section are not computed with EigenWave. 
Instead eigenvalues and eigenvectors of~\eqref{eq:discreteEigenvalueProblem} are computed using the Krylov-Schur algorithm in SLEPc~\cite{SLEPc2005}.

The computations presented in this section all take $L=\Delta$ in~\eqref{eq:helmholtzBVP} and~\eqref{eq:waveEquation} and use a time harmonic Gaussian source term having the form 
\ba
   f(\xv,t) = \alphag \cos(\omega t) \exp\bigl( -\betag^2  \| \xv-\xv_0\|^2 \bigr), 
   \label{eq:gaussianSource}
\ea
where $\alphag$ is the amplitude, $\xv_0=(x_0,y_0,z_0)$ denotes the center of the Gaussian, and the exponent coefficient $\betag$ 
is determines the approximate width of the Gaussian.
In the results that follow we choose $\alphag=\omega^2$, and $\betag=\omega$; this usually results in a solution with a 
maximum value that is roughly of size $1$.

\subsection{Experiments in a disk}

{
\newcommand{\drawContour}[7]{%
\begin{scope}[#1]
\draw(0.0,0) node[anchor=south west,xshift=-4pt,yshift=+0pt] {\trimfiga{fig/#2}{\figWidtha}};
  \draw(.5,.5) node[draw,fill=white,anchor=west,xshift=2pt,yshift=1pt] {\scriptsize #3};
\begin{scope}[xshift=-.2cm,yshift=+0pt]
  \draw (\xcb,\ycb) node[anchor=south west,xshift=0.25cm,yshift=.5cm,rotate=-90] {\trimfigcb{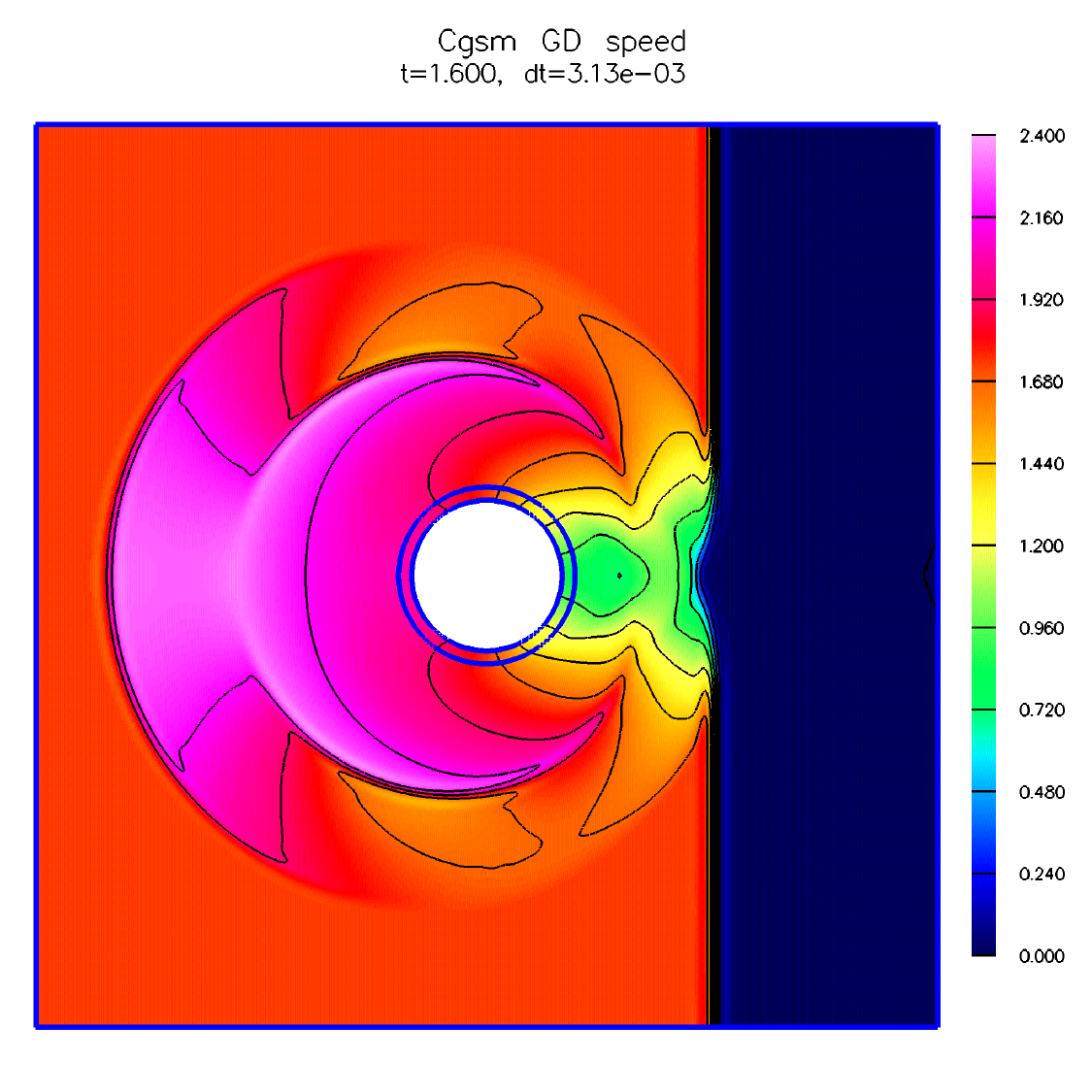}{\cbWidth}{\cbHeight}};
  \draw (.8,0) node[anchor=north,xshift=+3pt,yshift=+2pt] {\scriptsize $#6$};
  \draw (4.8,0) node[anchor=north,xshift=+0pt,yshift=+2pt] {\scriptsize $#7$};
\end{scope}
\end{scope}
}
\newcommand{\cbWidth}{.2cm}
\newcommand{\cbHeight}{4cm}
\newcommand{\xcb}{.5cm}
\newcommand{\ycb}{-.2cm}
\setlength{\ycbTop}{\ycb+\cbHeight}
\setlength{\ycbMid}{\ycb+\cbHeight*\real{.5}}
\newcommand{\trimfigcb}[3]{\includegraphics[width=#2, height=#3, clip, trim=17cm 2.35cm 1.65cm 2.35cm]{#1}}
\newcommand{\figWidtha}{4.5cm}
\newcommand{\trimfiga}[2]{\trimw{#1}{#2}{.11}{.115}{.11}{.11}}
\begin{figure}[htb]
\begin{center}
\resizebox{13.5cm}{!}{
\begin{tikzpicture}
   \useasboundingbox (0,.45) rectangle (15,4.5);  


   \begin{scope}[yshift=0cm]
     \figByWidth{   0}{-.2}{fig/sicGridG2}{5cm}[0.1][0.1][0.1][0.1]
     \drawContour{xshift=5.cm,yshift=0.00cm}{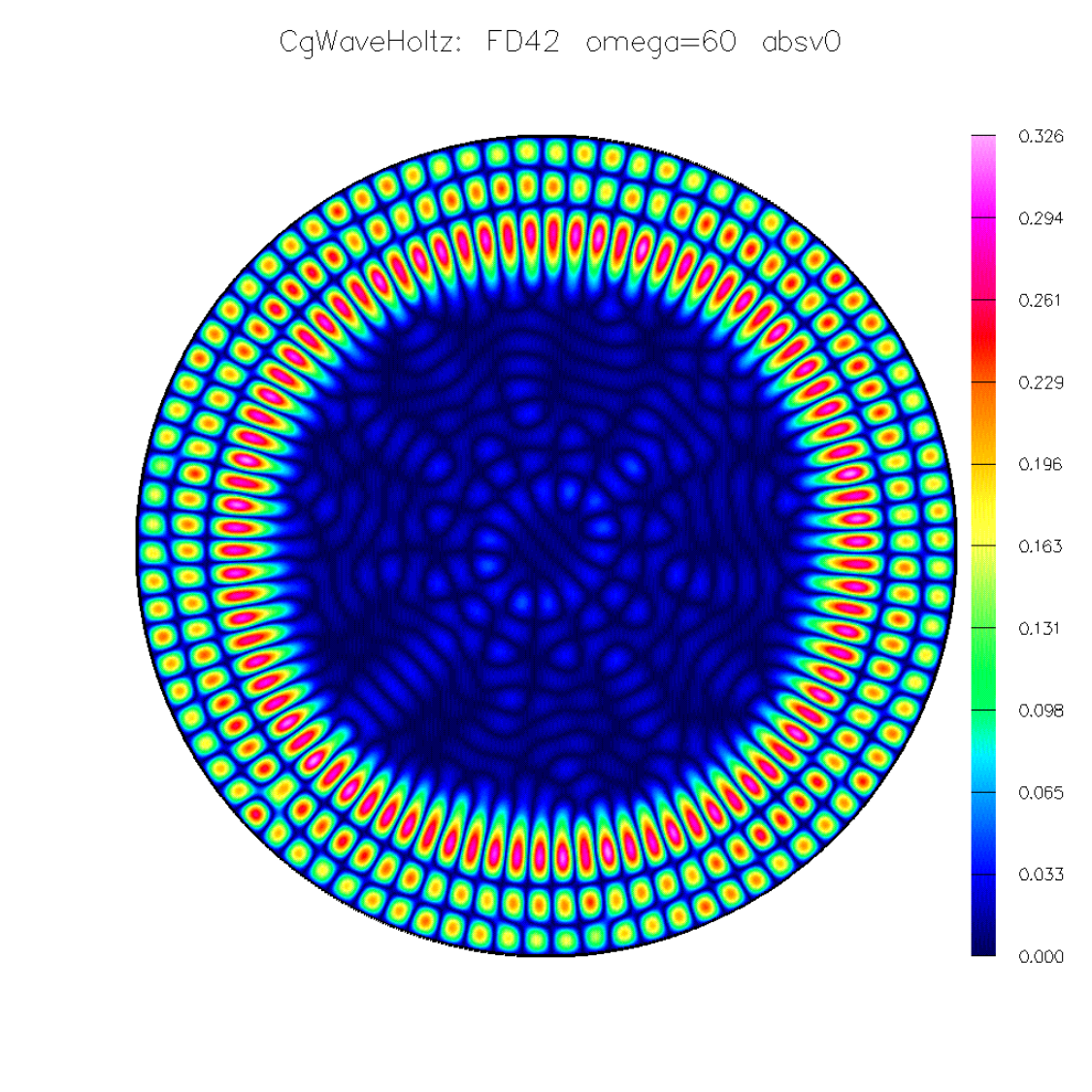}{$|v|$, $\omega=60$}{$v$}{$t=0.3$}{$-1.1$}{$0.33$};  
     \drawContour{xshift=10cm,yshift=0.00cm}{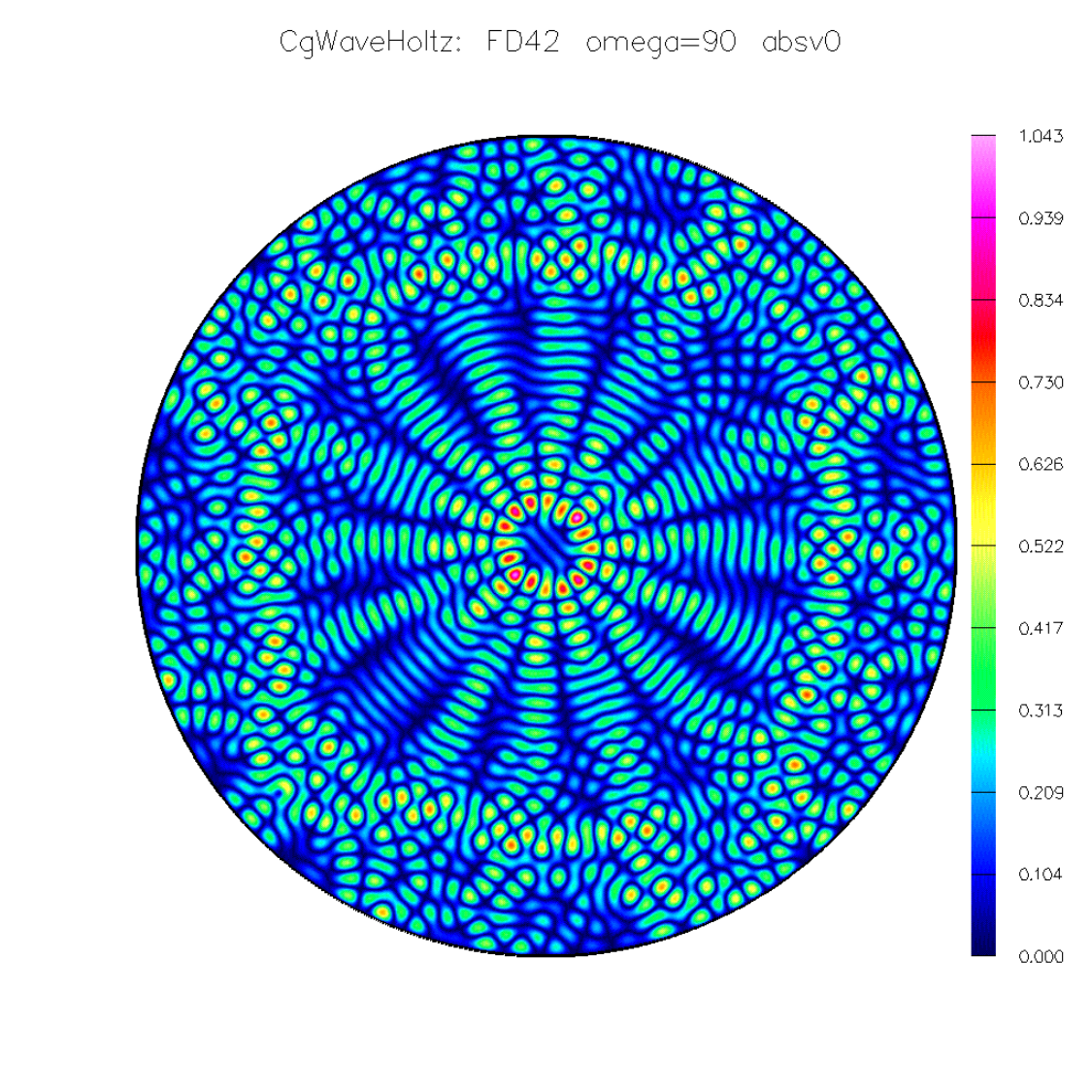}{$|v|$, $\omega=90$}{$v$}{$t=0.3$}{$0$}{$1.0$};  
        
   \end{scope}

\end{tikzpicture}
}
\end{center}
\caption{Left: coarse overset grid $\Gcd^{(2)}$ for the disk. 
  Middle and right: computed WaveHoltz solutions on grid $\Gcd^{(32)}$.
    }
\label{fig:diskGridAndSolutionFig}
\end{figure}
}

Helmholtz solutions are computed for a circular disk domain to demonstrate the
use of the WaveHoltz~scheme with an overset grid in two dimensions. 
The overset grid for the disk of radius $R=1$, as shown in the left plot of Figure~\ref{fig:diskGridAndSolutionFig},
consists of an annular boundary-fitted grid and a background Cartesian grid. 
Let $\Gcd^{(j)}$ denote the disk grid with target grid spacing $\ds^{(j)}=1/(10 j)$.
The middle and right plots of Figure~\ref{fig:diskGridAndSolutionFig} show sample WaveHoltz solutions
computed using deflation.


{
\newcommand{\drawContour}[7]{%
\begin{scope}[#1]
\draw(0.0,0) node[anchor=south west,xshift=-4pt,yshift=+0pt] {\trimfiga{fig/#2}{\figWidtha}};
\draw(3.8,.2) node[draw,fill=white,anchor=west,xshift=2pt,yshift=2pt,inner sep=3pt] {\scriptsize #3};
\begin{scope}[xshift=-12pt,yshift=-5pt]
  \draw (\xcb,\ycb) node[anchor=south west,xshift=0.25cm,yshift=.5cm,rotate=-90] {\trimfigcb{fig/colourBarLines}{\cbWidth}{\cbHeight}};
  \draw (.8,0) node[anchor=north,xshift=+3pt,yshift=+2pt] {\scriptsize $#6$};
  \draw (4.8,0) node[anchor=north,xshift=+0pt,yshift=+2pt] {\scriptsize $#7$};
\end{scope}
\end{scope}
}
\newcommand{\cbWidth}{.2cm}
\newcommand{\cbHeight}{4cm}
\newcommand{\xcb}{.5cm}
\newcommand{\ycb}{-.2cm}
\setlength{\ycbTop}{\ycb+\cbHeight}
\setlength{\ycbMid}{\ycb+\cbHeight*\real{.5}}
\newcommand{\trimfigcb}[3]{\includegraphics[width=#2, height=#3, clip, trim=17cm 2.35cm 1.65cm 2.35cm]{#1}}
\newcommand{\figWidtha}{4.2cm}
\newcommand{\trimfiga}[2]{\trimw{#1}{#2}{.12}{.12}{.12}{.12}}

\newcommand{\figSize}{5.5cm}

\begin{figure}[htb]

\begin{center}
\begin{tikzpicture}
   \useasboundingbox (0,0.25) rectangle (16,4.5);  

  \begin{scope}[yshift=0cm]  
    \drawContour{xshift=-.25cm,yshift=0.00cm}{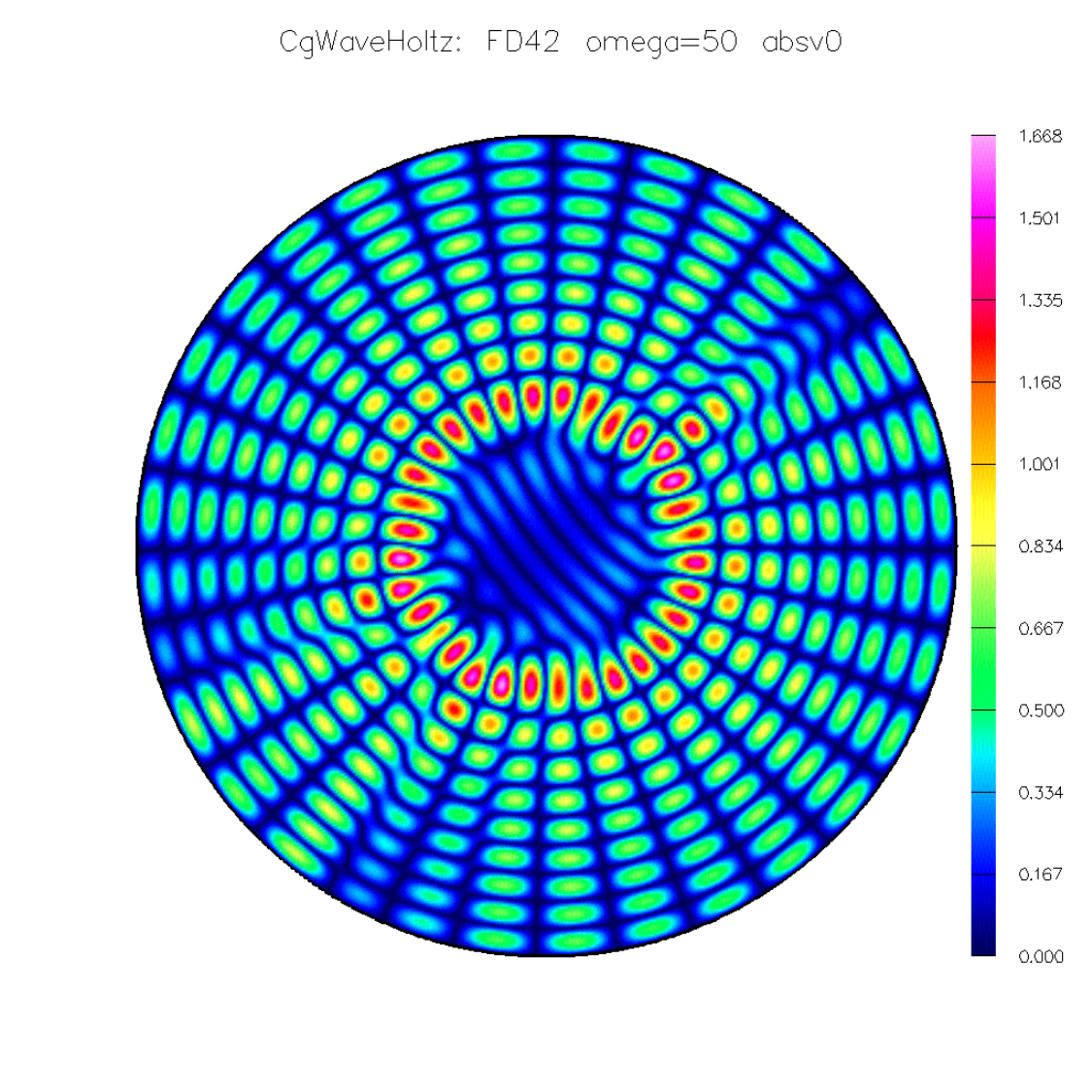}{$|v|$}{$v$}{$t=0.3$}{$0$}{$1.67$}    
   \end{scope}  
   \begin{scope}[xshift=4.75cm,yshift=-.25cm]

   \figByWidth{0.00}{0}{fig/diskG16O4Freq50Deflate64}{\figSize}[0][0][0][0];
   \figByWidth{5.75}{0}{fig/diskG16O4Freq50Deflate64FixPointMuFunction}{\figSize}[0][0][0][0];

  \end{scope}      

\end{tikzpicture}
\end{center}
\caption{Disk. Left: computed WaveHoltz solution for $\omega=50$ on grid $\Gcd^{(16)}$ at order four.
Middle: Convergence rates of the WaveHoltz fixed-point iterations (with theory) together with GMRES, AGMRES, BICGSTAB, and ABICGSTAB,
using  $\Ndeflate=125$ deflated eigenvectors. 
Right: WaveHoltz filter function $\beta$ together with the values of $\beta$ evaluated at the discrete eigenvalues $\lambda_{h,m}$.
 The black vertical line on the right graph indicates the value of the adjusted frequency $\omegaTilde$ used to correct for time discretization errors. 
}
\label{fig:diskG16Freq50Convergence}
\end{figure}

}

Figure~\ref{fig:diskG16Freq50Convergence} shows WaveHoltz convergence results for $\omega=50$ on grid $\Gcd^{(16)}$ at order four
using $N_p=2$ periods per time interval.
The Gaussian source~\eqref{eq:gaussianSource} was located at $\xv_0=[-0.25,-0.2]$.
The PPW was $20$ and the rule-of-thumb estimate~\eqref{eq:ppwPollution} with a relative error of $\eps=10^{-2}$ is also $20$ so that the computation was considered to be reasonably resolved.
The convergence rates are shown for the WaveHoltz fixed-point iteration (FPI), together with GMRES, AGMRES, BICGSTAB, and ABICGSTAB. A total of $\Na=125$ eigenvectors were deflated.
The FPI convergence is seen to be a good match to the theoretical asymptotic convergence rate (ACR) of $0.95$.
GMRES and AGMRES show almost identical convergence rates with an estimated convergence rate (CR) of about $0.66$.
BICGSTAB and ABICGSTAB are also very similar with an estimated CR of about $0.75$.
The right graph of Figure~\ref{fig:diskG16Freq50Convergence} shows WaveHoltz filter function $\beta$ together 
with the values of $\beta$ evaluated at the discrete eigenvalues $\lambda_{h,m}$.
The deflated eigenvalues are circled in black. 
From the distribution of eigenvalues it can be seen that almost any large value of $\omega$ will be relatively close to an eigenvalue.

It was previously noted that the discrete eigenvectors for an overset grid are accurate to $O(h^p)$ truncation errors but
not orthogonal to machine precision. This means there is some error in the projection steps in the direct eigenvector deflation (DEVD)
in Algorithm~\ref{alg:waveHoltzDirectDeflation} (really the left eigenvectors should also be used).
The question therefore arises as to the accuracy of the results using the DEVD method compared to
using the augmented-Krylov-eigenvector-deflation approach (AUKED).
For this problem the relative max-norm error between the WaveHoltz solution computed using DEVD (GMRES or BICGSTAB) and the solution of the
discretized Helmholtz problem solved with a sparse direct solver was $1.3\times 10^{-4}$. 
We conclude that DEVD is giving a good result given the PPW used.
The corresponding relative max-norm error using AUKED (AGMRES or ABICGSTAB) was $2.4\times10^{-8}$. 
In exact arithmetic the AUKED answer should exactly converge to the solution of the discretized Helmholtz problem.
The reason this is not so is due to finite precision and the conditioning of the linear systems.
Note that the WaveHoltz linear system in~\eqref{eq:WaveHoltzMatrixEquation} involves a different matrix from that formed from the
discretized Helmholtz problem.

{
\newcommand{\figSize}{6cm}

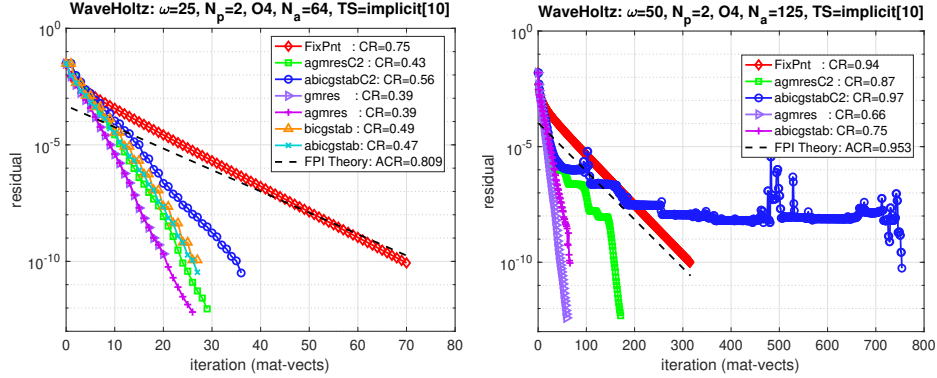
\begin{figure}[htb]

\begin{center}
\begin{tikzpicture}
   \useasboundingbox (0,.65) rectangle (13,4.95);  

  \begin{scope}
  \figByWidth{0.25}{0}{fig/diskG16O4Freq25Deflate32EVC}{\figSize}[0][0][0][0];
  \figByWidth{6.5}{0}{fig/diskG16O4Freq50Deflate64EVC}{\figSize}[0][0][0][0];
  \end{scope}      

\end{tikzpicture}
\end{center}
\caption{Disk: Comparing deflation with fine-grid eigenvectors to deflation with coarse-grid eigenvectors, C2 = coarse grid eigenvectors.
Left: $\omega=25$ (coarse grid eigenvectors are well resolved).
Right: $\omega=50$ (coarse grid eigenvectors are not well resolved).
Startup iterations for C2 computations are not included.
}
\label{fig:disk16Freq50ConvergenceC2}
\end{figure}
}

The use of coarse grid eigenvectors for deflation is demonstrated in Figure~\ref{fig:disk16Freq50ConvergenceC2}.
The fine grid is $\Gcd^{(16)}$ and the coarse grid is $\Gcd^{(8)}$.
Results using coarse grid eigenvectors are labeled agmresC2 and abicgstabC2 (C2 means coarsened by a factor of $2$).
The left graph of Figure~\ref{fig:disk16Freq50ConvergenceC2} shows results for $\omega=25$ when the 
the coarse grid eigenvectors are well resolved (the computed PPW is 20 on the coarse grid with rule-of-thumb estimate~\eqref{eq:ppwPollution}  of $\PPW_4=17$).
AGMRESC2 shows good convergence with $\CR\approx 0.43$ compared to GMRES with $\CR\approx 0.39$.
ABICGSTABC2 also shows good convergence with $\CR\approx 0.56$ compared to BICGSTAB with $\CR\approx 0.47$.
The right graph of Figure~\ref{fig:disk16Freq50ConvergenceC2} shows results for $\omega=50$ when
the coarse grid eigenvectors are not well resolved (the computed PPW is 10 on the coarse grid with rule-of-thumb estimate of $\PPW_4=20$).
In this case the augmented Krylov schemes with coarse grid eigenvectors do not converge very well compared to the use of fine grid
eigenvectors. In particular ABICGSTABC2 shows difficulty in converging past a tolerance of about $10^{-8}$.

{
\newcommand{\figSize}{6cm}

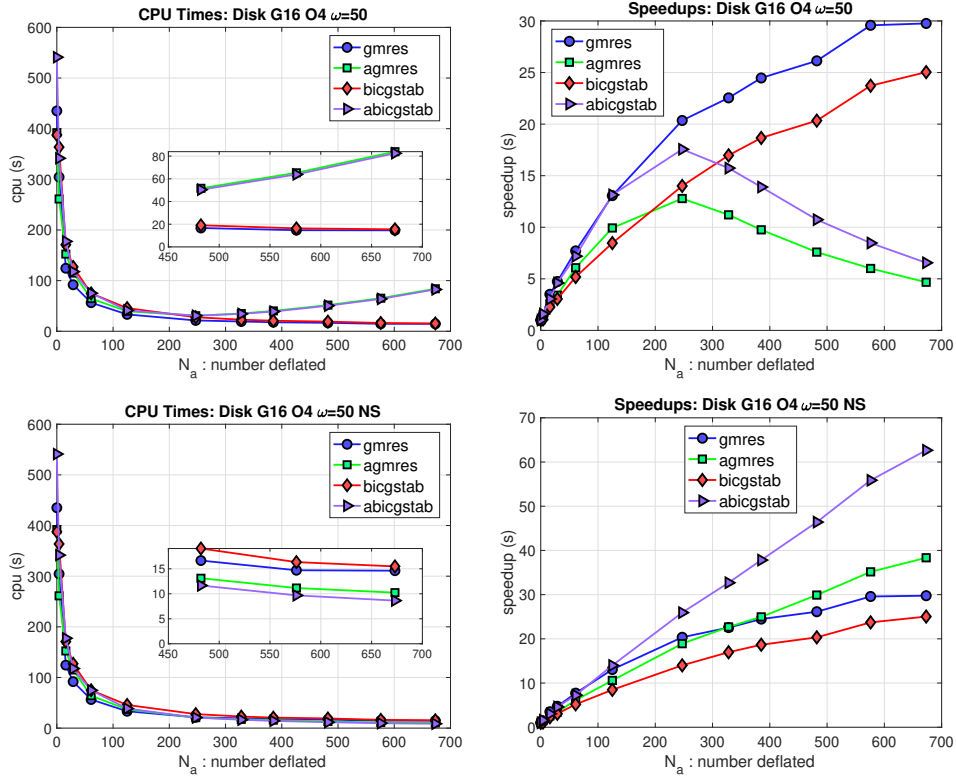
\begin{figure}[htb]

\begin{center}
\begin{tikzpicture}
   \useasboundingbox (0,.5) rectangle (13,10.5);  

 \begin{scope}[yshift=5.25cm]
    \figByWidth{0.00}{0}{fig/diskG16O4Freq50CpuVersusDeflated}{\figSize}[0][0][0][0]
    \figByWidth{6.5}{0}{fig/diskG16O4Freq50SpeedupVersusDeflated}{\figSize}[0][0][0][0]
  \end{scope}    
  \begin{scope}
    \figByWidth{0.00}{0}{fig/diskG16O4Freq50NSCpuVersusDeflated}{\figSize}[0][0][0][0]
    \figByWidth{6.5}{0}{fig/diskG16O4Freq50NSSpeedupVersusDeflated}{\figSize}[0][0][0][0]
  \end{scope}      

\end{tikzpicture}
\end{center}
\caption{Disk CPU times and speedups. Top: including startup times. Bottom: excluding startup times (NS=no startup).
}
\label{fig:disk16Freq50Speedups}
\end{figure}
}

CPU time speedups that are obtained using eigenvector deflation are indicated in Figure~\ref{fig:disk16Freq50Speedups} for computations
on the disk grid $\Gcd^{(16)}$ with $\omega=50$.
The left column shows CPU times in seconds for the four schemes GMRES, AGMRES, BICGSTAB, and ABICGSTAB to reach
a convergence tolerance of $10^{-10}$. The right column shows the CPU speedup of each scheme compared to the CPU time 
without deflation, $\Na=0$.
The top row shows results that include the startup costs for the augmented routines; this being the cost of the QR factorization
in algorithms~\ref{alg:augmentedGmresWithGivens} and~\ref{alg:augmentedBiCGStab}.
The bottom row shows results not including the startup times (which can be precomputed, independent of the forcing $f$).
When not including startup times, all schemes show good speedups out to $\Na \approx 700$.
Augmented BICGSTAB tends to perform relatively poorly for small values of $\Na$ but behaves well for larger $\Na$. 
With startup costs included the augmented routines reach a maximum speedup at $\Na\approx 250$.
For small $\Na$, when the number of iterations is large, the primary cost of memory for GMRES and AGMRES is storing the Krylov space vectors.
For large $\Na$, when the number of iterations is small, the primary cost in memory is for storing the augmented vectors. In both cases this memory cost could be alleviated by compression techniques as described for the single grid case.

\newcommand{\Gcde}{\Gc_{\rm de}}
\subsection{Penrose unilluminable room}

As a next example, we solve the Helmholtz problem for the Penrose unilluminable room~\cite{Fukushima2015LightPI}.
The geometry, shown in Figure~\ref{fig:doubleEllipseGridFig}, is designed so that the 
some of the alcoves, two at the top and two at the bottom of the domain, remain dark (or quiet) when there is a light source (or sound source) in the interior.
The design is based on two ellipses of different sizes. 
Two smaller half-ellipses, with semi-axes $(a_1,b_1)=(2,1)$, are located at the top and bottom. 
Two larger half-ellipses, with semi-axes $(a_2,b_2)=(3,6)$, are placed on the left and right.
The left and right ends of the smaller ellipses are located at the foci of the larger ellipses.

{
\newcommand{\drawContour}[7]{%
\begin{scope}[#1]
\draw(0.0,0) node[anchor=south west,xshift=-4pt,yshift=+0pt] {\trimfiga{#2}{\figWidtha}};
\draw(.5,.5) node[draw,fill=white,anchor=west,xshift=2pt,yshift=1pt] {\scriptsize #3};
\begin{scope}[xshift=-3pt,yshift=0pt]
  \draw (\xcb,\ycb) node[anchor=south west,xshift=0.25cm,yshift=.5cm,rotate=-90] {\trimfigcb{fig/colourBarLines}{\cbWidth}{\cbHeight}};
  \draw (.8,0) node[anchor=north,xshift=+3pt,yshift=+2pt] {\scriptsize $#6$};
  \draw (4.8,0) node[anchor=north,xshift=+0pt,yshift=+2pt] {\scriptsize $#7$};
\end{scope}
\end{scope}
}
\newcommand{\cbWidth}{.2cm}
\newcommand{\cbHeight}{4cm}
\newcommand{\xcb}{.5cm}
\newcommand{\ycb}{-.2cm}
\setlength{\ycbTop}{\ycb+\cbHeight}
\setlength{\ycbMid}{\ycb+\cbHeight*\real{.5}}
\newcommand{\trimfigcb}[3]{\includegraphics[width=#2, height=#3, clip, trim=17cm 2.35cm 1.65cm 2.35cm]{#1}}
\newcommand{\figWidtha}{5.85cm}
\newcommand{\trimfiga}[2]{\trimh{#1}{#2}{.17}{.17}{.1}{.1}}
\newcommand{\figSize}{5.5cm}
\newcommand{\figh}{5.5cm}
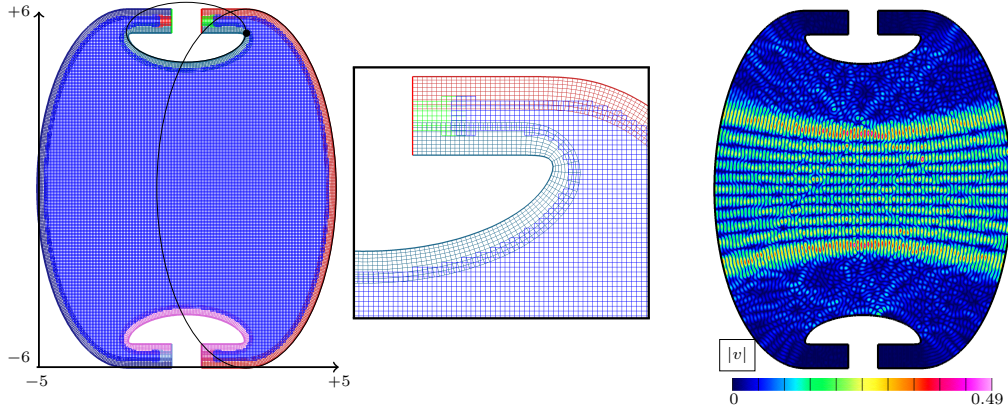
\begin{figure}[htb]
\begin{center}
\resizebox{13.cm}{!}{
\begin{tikzpicture}
  \useasboundingbox (0,.2) rectangle (15,5.7);  

  \figByWidth{0.0}{-6pt}{fig/darkCornerRoomGridG2}{4.7cm}[0.1][0.1][0.][0.]
  \draw[thick,black,->,yshift=0pt] (0,0) node[anchor=north] {\scriptsize $-5$} -- (4.65,0.00) node[anchor=north] {\scriptsize $+5$};  
  \draw[thick,black,->,xshift=1pt] (0,0) node[anchor=east,yshift=4pt]  {\scriptsize $-6$} -- (0.00,5.47) node[anchor=east ] {\scriptsize $+6$};  

  \figByWidthb{4.95}{.75}{fig/darkCornerRoomGridG2Zoom}{4.5cm}[0.][0.][0.][0.15]

  \begin{scope}[scale=0.46,xshift=7.02cm,yshift=5.97cm]
    \draw (0,0) ellipse (3cm and 6cm);    
    \draw[black,fill=black] (0,5.196cm) circle (3pt); 
  \end{scope}
  \begin{scope}[scale=0.46,xshift=5cm,yshift=11.2cm]
    \draw (0,0) ellipse (2cm and 1cm);    
  \end{scope}
  
  \begin{scope}[xshift=10cm,yshift=0]

     \drawContour{xshift=-2pt,yshift=-10pt}{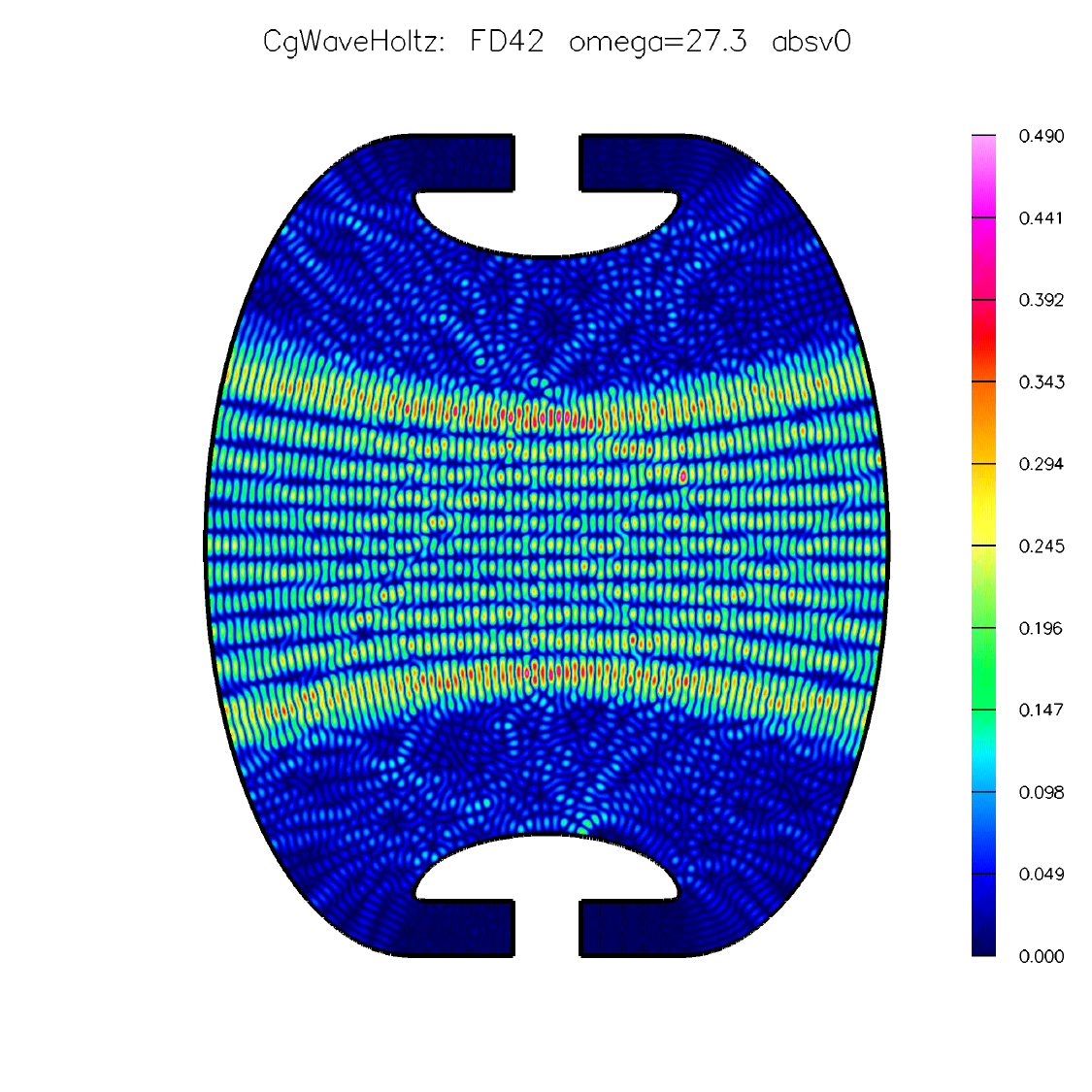}{$|v|$}{$v$}{$t=0.3$}{$0$}{$0.49$}
     

  \end{scope} 
\end{tikzpicture}
} 
\end{center}
\caption{Left: double ellipse geometry and overset grid $\Gcde^{(2)}$. Middle: closeup of a portion of the grid.
Right: computed Helmholtz solution for $\omega=27.3$ with a Gaussian source at $(2,1)$.
 }
\label{fig:doubleEllipseGridFig}
\end{figure}
}

The overset grid for the domain is shown in Figure~\ref{fig:doubleEllipseGridFig} (left and middle).
The grid, denoted by $\Gcde^{(j)}$ with target grid spacing $\ds^{(j)}=1/(10 j)$, consists of a total of nine component grids.  Four component grids are placed to fit the curved elliptical boundaries with four small Cartesian grids used to fit the straight portions of the boundaries in the alcoves (see middle image).  The ninth component grid is a large background Cartesian grid covering the bulk of the domain.
Figure~\ref{fig:doubleEllipseGridFig} (right) shows a sample solution computed by the WaveHoltz algorithm,
using a Gaussian source~\eqref{eq:gaussianSource} with  
$\omega=27.3$ and located at $\xv_0=(2,1)$.
The forcing excites an harmonic mode that is active primarily in the central region of the room, leaving the four alcoves quiet.


{
\newcommand{\drawContour}[7]{%
\begin{scope}[#1]
\draw(0.0,0) node[anchor=south west,xshift=-4pt,yshift=+0pt] {\trimfiga{#2}{\figWidtha}};
 \draw(3.7,.3) node[draw,fill=white,anchor=west,xshift=2pt,yshift=2pt,inner sep=3pt] {\footnotesize #3};
\begin{scope}[xshift=-12pt,yshift=-5pt]
  \draw (\xcb,\ycb) node[anchor=south west,xshift=0.25cm,yshift=.5cm,rotate=-90] {\trimfigcb{fig/colourBarLines}{\cbWidth}{\cbHeight}};
  \draw (.8,0) node[anchor=north,xshift=+3pt,yshift=+2pt] {\scriptsize $#6$};
  \draw (4.8,0) node[anchor=north,xshift=+0pt,yshift=+2pt] {\scriptsize $#7$};
\end{scope}
\end{scope}
}
\newcommand{\cbWidth}{.2cm}
\newcommand{\cbHeight}{4cm}
\newcommand{\xcb}{.5cm}
\newcommand{\ycb}{-.2cm}
\setlength{\ycbTop}{\ycb+\cbHeight}
\setlength{\ycbMid}{\ycb+\cbHeight*\real{.5}}
\newcommand{\trimfigcb}[3]{\includegraphics[width=#2, height=#3, clip, trim=17cm 2.35cm 1.65cm 2.35cm]{#1}}
\newcommand{\figWidtha}{4.cm}
\newcommand{\trimfiga}[2]{\trimw{#1}{#2}{.18}{.18}{.11}{.11}}

\newcommand{\figSize}{5.75cm}

\begin{figure}[htb]

\begin{center}
\begin{tikzpicture}
   \useasboundingbox (0,.5) rectangle (15.75,5.4);  

  \begin{scope}[yshift=0.4cm]  
    \drawContour{xshift=-.5cm,yshift=0.00cm}{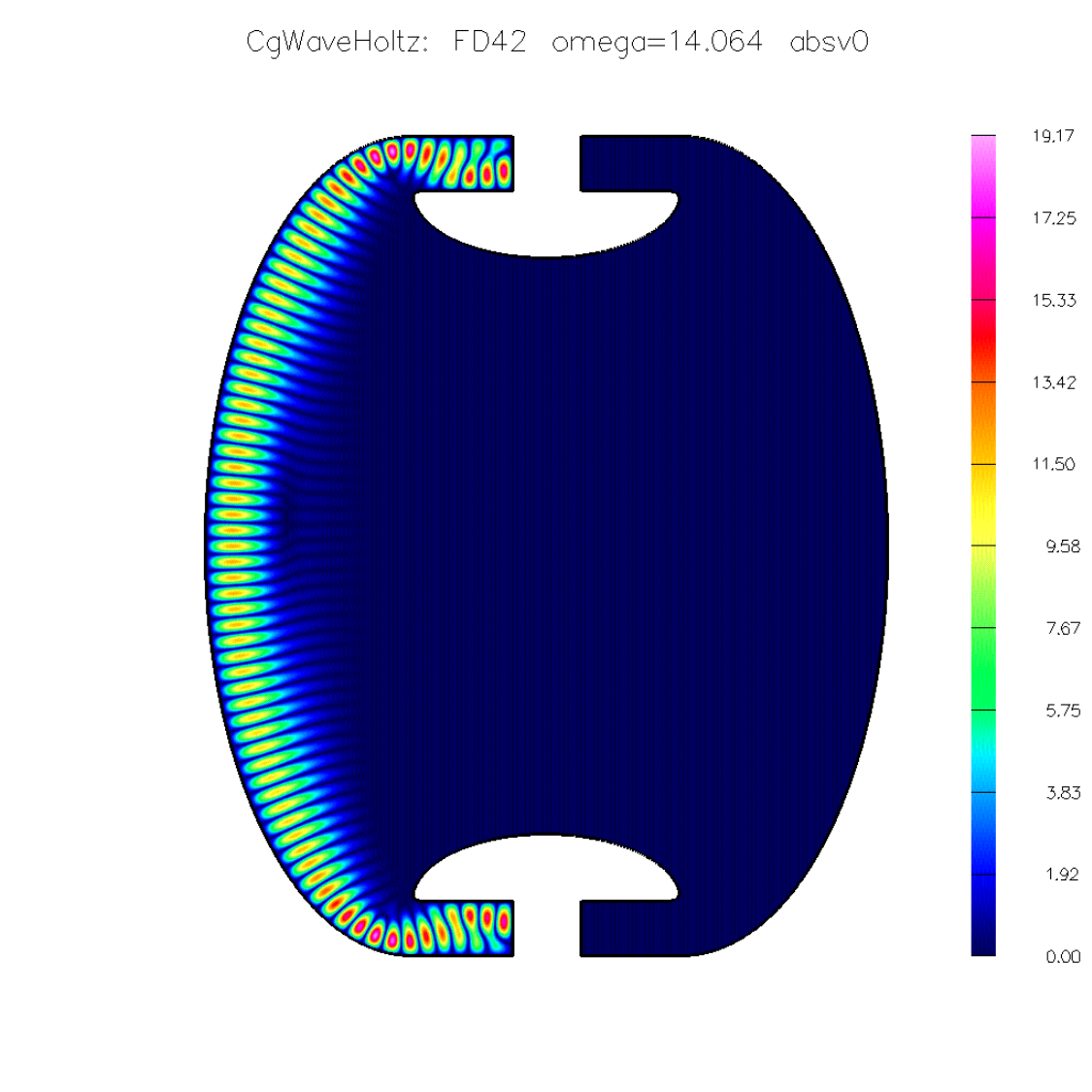}{$|v|$}{$v$}{$t=0.3$}{$0$}{$19$}    
  \end{scope}  
  \begin{scope}[xshift=4.05cm,yshift=0cm]
   \figByWidth{0.25}{0}{fig/darkCornerRoomG8O4Freq14}{\figSize}[0][0][0][0];
   \figByWidth{6.0}{0}{fig/darkCornerRoomG8O4Freq14FixPointMuFunction}{\figSize}[0][0][0][0];

  \end{scope}      

\end{tikzpicture}
\end{center}
\caption{Unliminable room. Left: WaveHoltz solution for $\omega=14.06$, computed on grid $\Gcde^{(8)}$ at order $4$.
Middle: Convergence rates of the WaveHoltz fixed-point iterations (with theory) together with GMRES, AGMRES, BICGSTAB, and ABICGSTAB,
using  $\Ndeflate=363$ deflated eigenvectors. 
Right: WaveHoltz filter function $\beta$ together with the values of $\beta$ evaluated at the discrete eigenvalues $\lambda_{h,m}$.
 The black vertical line on the right graph indicates the value of the adjusted frequency $\omegaTilde$ used to correct for time discretization errors. 
}
\label{fig:doubleEllipseFreq14}
\end{figure}

}

Figure~\ref{fig:doubleEllipseFreq14}
shows convergence results for GMRES, AGMRES, BICSTAB and ABICGSTAB for a case 
with  $\omega=14.06$ (a value close to resonance) and a Gaussian source located in the lower left alcove at $\xv_0=[-1.4,-5.6]$. 
The forcing in this case excites a surface mode that primarily resides along the left side of the room.
The solution is computed on grid $\Gcde^{(8)}$ to fourth order accuracy
using implicit time-stepping with $10$ time-steps per period and $N_p=2$.
A total of $363$ eigenvectors were used for deflation.
For this problem the solution was computed using $\PPW=36$ 
points-per-wavelength, while the rule of thumb estimate~\eqref{eq:ppwPollution} ($\eps=10^{-2}$) was $\PPW_4=23$; this suggests the computation should be reasonably resolved.
The convergence of the FPI is approaching the theoretical ACR.
The GMRES and AGMRES results show similar convergence rates with $\CR\approx 0.62$. The BICGSTAB and ABICGSTAB both have
convergence rates of about $\CR\approx 0.73$.

\newcommand{\Gcke}{\Gc_{\rm ke}}
\subsection{Sixteen knife edges}


{
\newcommand{\labelSize}{\normalsize}
\newcommand{\drawContour}[7]{%
\begin{scope}[#1]
\draw(0.0,0) node[anchor=south west,xshift=-4pt,yshift=+0pt] {\trimfiga{#2}{\figWidtha}};
\draw(.5,.5) node[draw,fill=white,anchor=west,xshift=2pt,yshift=1pt,inner sep=2pt] {\labelSize #3};
\begin{scope}[xshift=0.5cm,yshift=-2pt]
  \draw (\xcb,\ycb) node[anchor=south west,xshift=0.25cm,yshift=.5cm,rotate=-90] {\trimfigcb{fig/colourBarLines}{\cbWidth}{\cbHeight}};
  \draw (.8,0) node[anchor=north,xshift=+3pt,yshift=+2pt] {\labelSize $#6$};
  \draw (6.8,0) node[anchor=north,xshift=+0pt,yshift=+2pt] {\labelSize $#7$};
\end{scope}
\end{scope}
}
\newcommand{\cbWidth}{.2cm}
\newcommand{\cbHeight}{6cm}
\newcommand{\xcb}{.5cm}
\newcommand{\ycb}{-.2cm}
\setlength{\ycbTop}{\ycb+\cbHeight}
\setlength{\ycbMid}{\ycb+\cbHeight*\real{.5}}
\newcommand{\trimfigcb}[3]{\includegraphics[width=#2, height=#3, clip, trim=17cm 2.35cm 1.65cm 2.35cm]{#1}}
\newcommand{\figWidtha}{8.27cm}
\newcommand{\trimfiga}[2]{\trimw{#1}{#2}{.12}{.12}{.28}{.28}}

\begin{figure}[htb]
\begin{center}
\resizebox{12.cm}{!}{
\begin{tikzpicture}
   \useasboundingbox (0,.2) rectangle (18,9.8);  

   \begin{scope}[xshift=0cm,yshift=5.25cm]
	   \figByHeight{0.04}{-.001}{fig/sixteenKnivesG64}{4.5cm}[0.03][0.03][0.24][0.24]
	   \figByHeightb{8.75}{0}{fig/sixteenKnivesG64Zoom1}{4.45cm}[0.][0.][0.][0.]
	   \figByHeightb{13.75}{0}{fig/sixteenKnivesG64Zoom2}{4.45cm}[0.][0.][0.][0.]

	  \draw[thick,black,->,yshift=0pt] (0,0) -- (8.4,0.00);  
	  \draw[thick,black,->,xshift=-1pt] (0,0) -- (0.00,4.8);     

	   \draw[thick,black,-,xshift=0pt] (0.00,.1) -- (0.00,-.1) node[anchor=north] {\labelSize $-.5$};
	   \draw[thick,black,-,xshift=0pt] (4.07,.1) -- (4.07,-.1) node[anchor=north] {\labelSize $0$};
	   \draw[thick,black,-,xshift=0pt] (8.10,.1) -- (8.10,-.1) node[anchor=north] {\labelSize $0.5$};

	   \draw[thick,black,-,xshift=0pt] (.1,0.0) -- (-.1,0.0) node[anchor=east] {\labelSize $0$};
	   \draw[thick,black,-,xshift=0pt] (.1,4.5) -- (-.1,4.5) node[anchor=east] {\labelSize $.55$};
   \end{scope}

   \begin{scope}[xshift=.25cm,yshift=-.29cm]
     \drawContour{xshift= 0cm,yshift=0.cm}{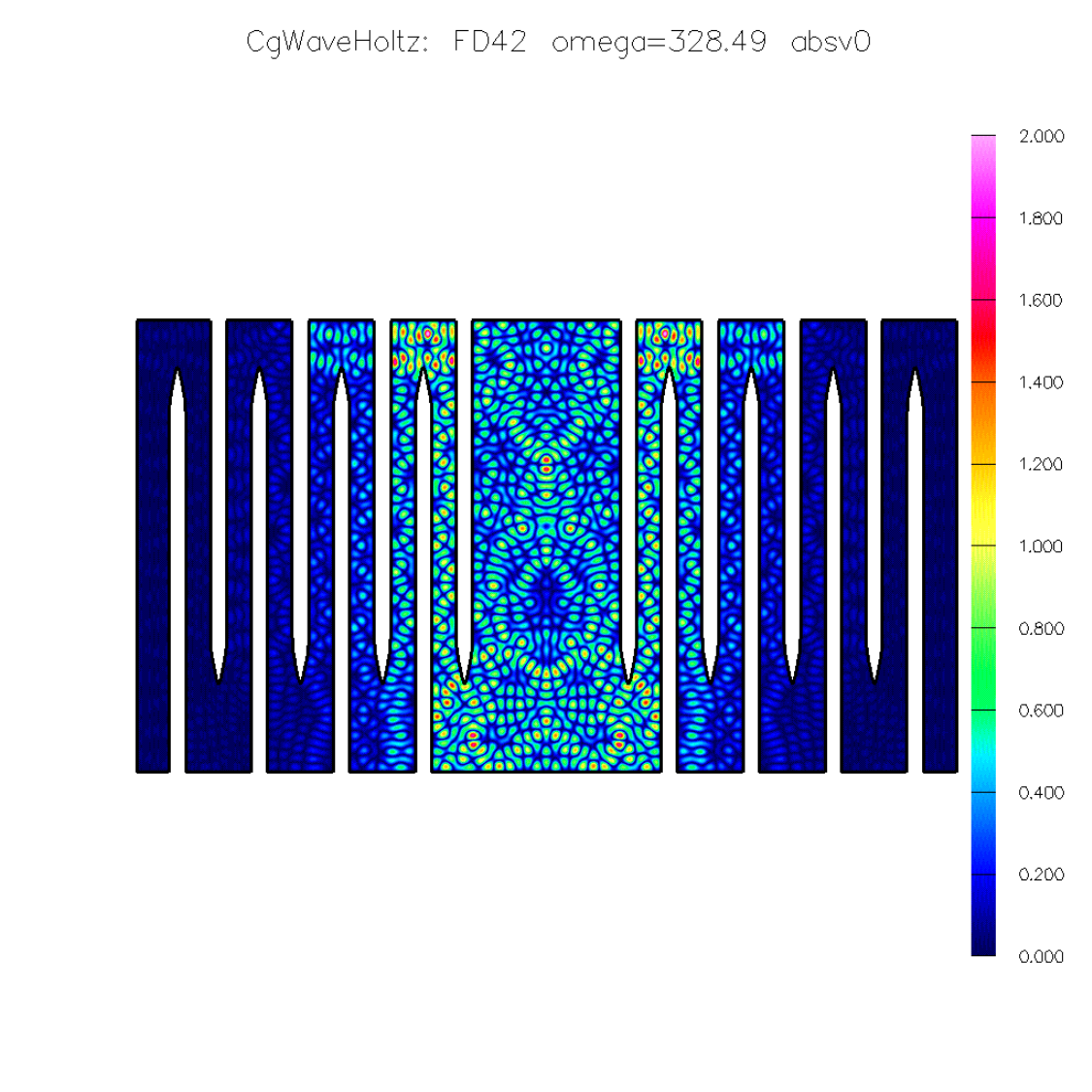}{$|v|$}{$v$}{$t=1.0$}{$0.$}{$2.0$};     
     \drawContour{xshift=8.5cm,yshift=0.cm}{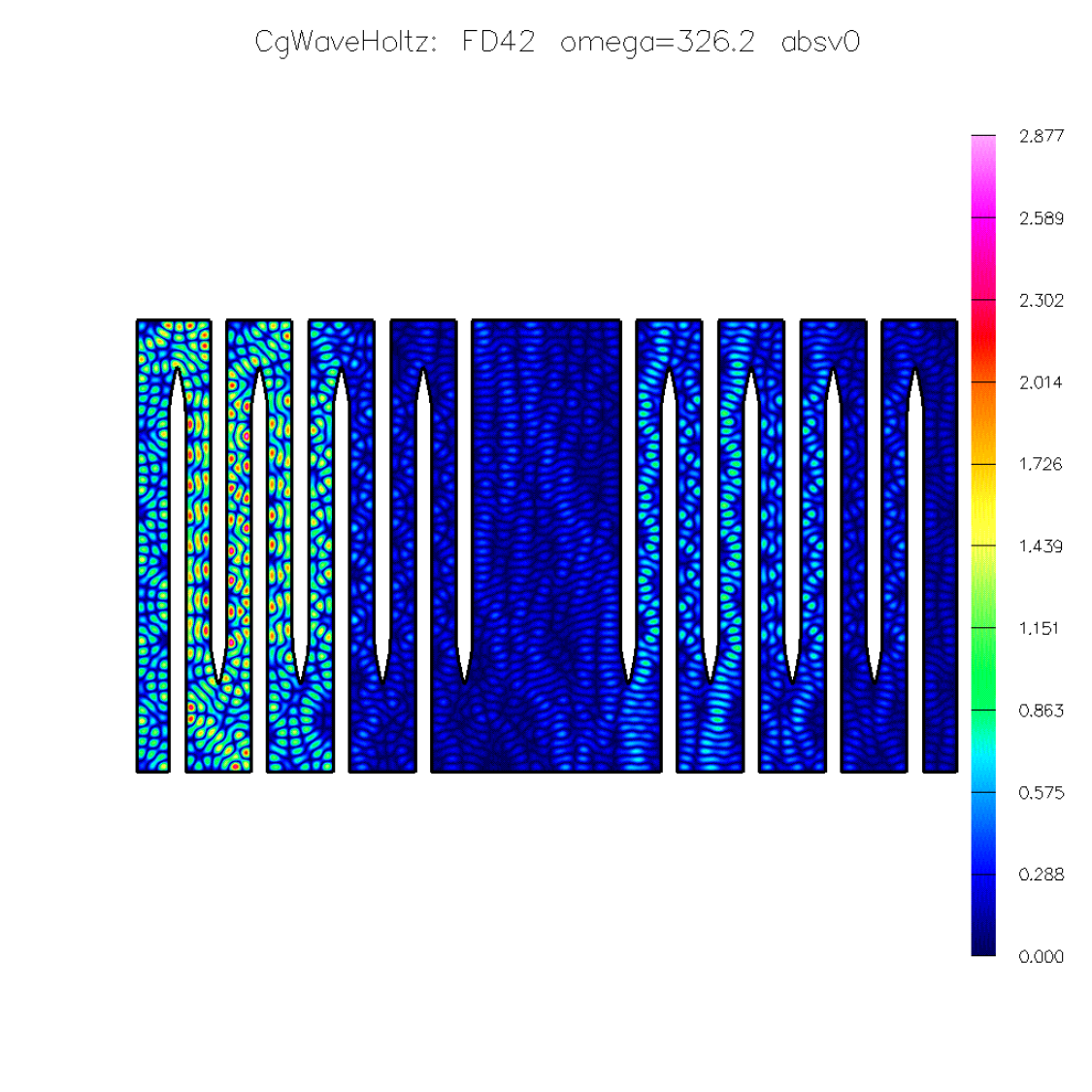}{$|v|$}{$v$}{$t=1.0$}{$0.$}{$2.88$};    
   \end{scope}   

\end{tikzpicture}
}
\end{center}
\caption{
    Sixteen knife edges. Top: overset grid $\Gcke^{(64)}$ and magnified views. 
    Bottom: WaveHoltz solutions computed on grid $\Gcke^{(128)}$, order four.
    Bottom left: $\omega=328.49$, source in the middle $\xv_0=(0,0.3)$.
    Bottom right: $\omega=326.2$, source in the upper left $\xv_0=(-.47,0.52)$.
    }
\label{fig:sixteenKnivesGridAndContours}
\end{figure}
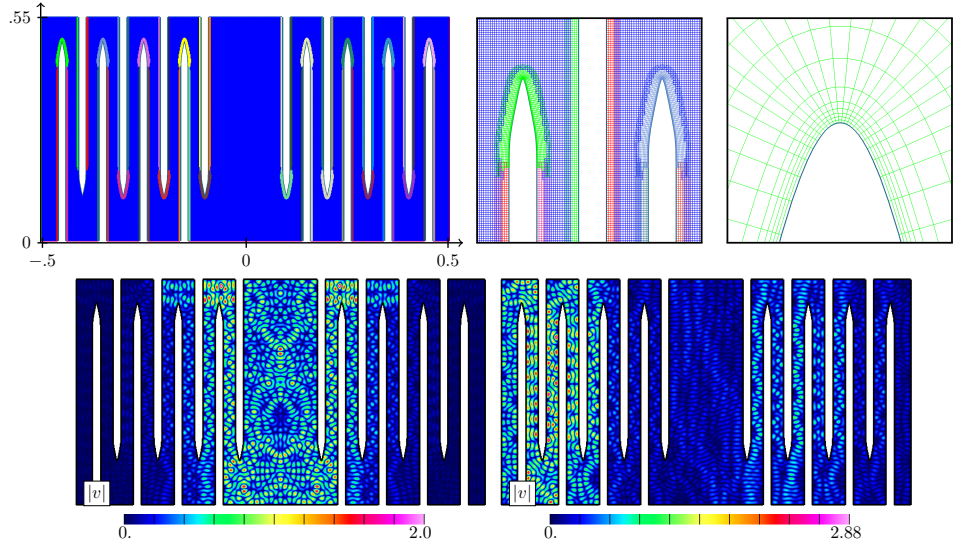

}

Helmholtz solutions are computed for a complex geometry consisting of sixteen knife edges as shown
in Figure~\ref{fig:sixteenKnivesGridAndContours}. This is an example where implicit time stepping is especially useful 
since the grid contains small cells to resolve the tips of the knife edges. A standard explicit time-stepping scheme
would require a small time-step everywhere due to the presence of a relatively few tiny cells.


{
\newcommand{\drawContour}[7]{%
\begin{scope}[#1]
\draw(0.0,0) node[anchor=south west,xshift=-4pt,yshift=+0pt] {\trimfiga{#2}{\figWidtha}};
\begin{scope}[xshift=0cm,yshift=-5pt]
  \draw (\xcb,\ycb) node[anchor=south west,xshift=0.25cm,yshift=.5cm,rotate=-90] {\trimfigcb{fig/colourBarLines}{\cbWidth}{\cbHeight}};
  \draw (.8,0) node[anchor=north,xshift=+3pt,yshift=+2pt] {\scriptsize $#6$};
  \draw (4.8,0) node[anchor=north,xshift=+0pt,yshift=+2pt] {\scriptsize $#7$};
\end{scope}
\end{scope}
}
\newcommand{\cbWidth}{.2cm}
\newcommand{\cbHeight}{4cm}
\newcommand{\xcb}{.5cm}
\newcommand{\ycb}{-.2cm}
\setlength{\ycbTop}{\ycb+\cbHeight}
\setlength{\ycbMid}{\ycb+\cbHeight*\real{.5}}
\newcommand{\trimfigcb}[3]{\includegraphics[width=#2, height=#3, clip, trim=17cm 2.35cm 1.65cm 2.35cm]{#1}}
\newcommand{\figWidtha}{5cm}
\newcommand{\trimfiga}[2]{\trimw{#1}{#2}{.11}{.12}{.285}{.285}}

\newcommand{\figSize}{5.25cm}

\begin{figure}[htb]

\begin{center}
\begin{tikzpicture}
   \useasboundingbox (0,.5) rectangle (16,4.5);  

  \begin{scope}[yshift=.8cm]  
    \drawContour{xshift=-.5cm,yshift=0.00cm}{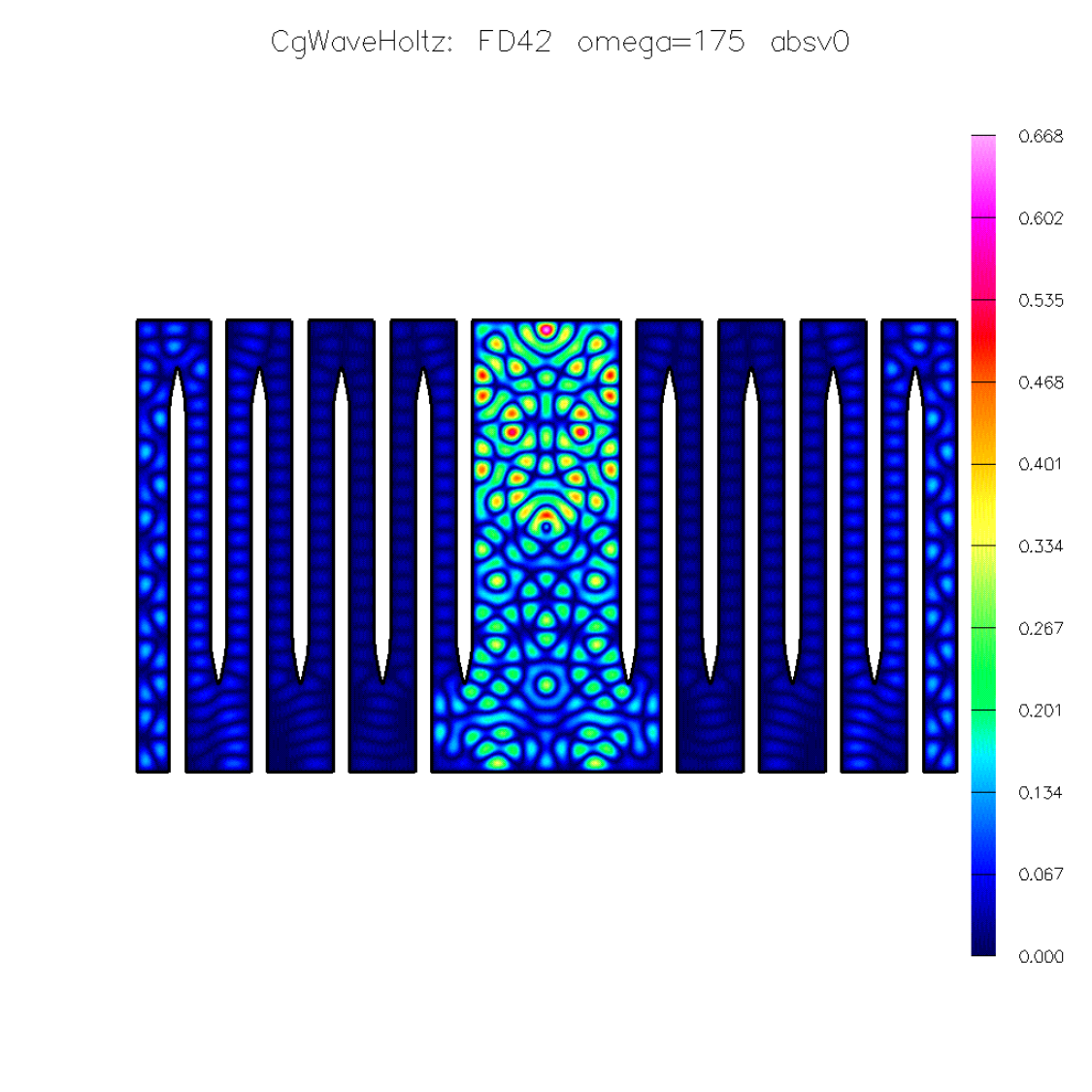}{$|v|$}{$v$}{$t=0.3$}{$0$}{$0.67$}    
   \end{scope}  
   \begin{scope}[xshift=4.75cm,yshift=0cm]
   \figByWidth{0.25}{0}{fig/sixteenKnivesGride64O4Freq175}{\figSize}[0][0][0][0];
   \figByWidth{5.75}{0}{fig/sixteenKnivesGride64O4Freq175FixPointMuFunction}{\figSize}[0][0][0][0];

  \end{scope}      

\end{tikzpicture}
\end{center}
\caption{Sixteen knives. Left: WaveHoltz solution for $\omega=175$, computed on grid $\Gcke^{(64)}$ at order $4$.
Middle: Convergence rates of the WaveHoltz fixed-point iterations (with theory) together with GMRES, AGMRES, BICGSTAB, and ABICGSTAB,
using  $\Ndeflate=216$ deflated eigenvectors. 
Right: WaveHoltz filter function $\beta$ together with the values of $\beta$ evaluated at the discrete eigenvalues $\lambda_{h,m}$.
 The black vertical line on the right graph indicates the value of the adjusted frequency $\omegaTilde$ used to correct for time discretization errors. 
}
\label{fig:sixteenKnivesFreq175}
\end{figure}

}

The overset grid for the geometry is shown in the top plots of Figure~\ref{fig:sixteenKnivesGridAndContours}.
Cartesian grids cover the sides of the knife edges while curvilinear grids are used around the tips.
These grids are embedded in a background Cartesian grid. 
The overset grid of resolution factor $j$, denoted by  $\Gcke^{(j)}$, and has target grid spacing $\ds=1/(10 j)$.
The bottom plots of Figure~\ref{fig:sixteenKnivesGridAndContours}
show sample solutions computed using WaveHoltz and deflation.

Figure~\ref{fig:sixteenKnivesFreq175} shows convergence results for GMRES, AGMRES, BICSTAB and ABICGSTAB for a case 
with  $\omega=175$ and a Gaussian source located at $\xv_0=[0,0.3]$. The solution is computed on grid $\Gcke^{(64)}$ to fourth order accuracy
using implicit time-stepping with $10$ time-steps per period and $N_p=2$.
A total of $216$ eigenvectors were used for deflation.
For this problem the actual and rule of thumb estimates ($\eps=10^{-2}$) for the points-per-wavelength were both $\PPW=23$; this suggests the computation should be reasonably resolved.
The convergence of the FPI is approaching the theoretical ACR.
The GMRES and AGMRES results show similar convergence rates with $\CR\approx 0.57$. The BICGSTAB and ABICGSTAB both have
convergence rates of about $\CR\approx 0.68$.

\section{Conclusions \label{sec:conclusion}}
We have carried out a systematic numerical study of eigenvector deflation for accelerating the WaveHoltz method. Both the direct deflation approach (DEVD) and the augmented-Krylov approach (AUKED) markedly reduce the iteration count relative to undeflated WaveHoltz and standard Krylov acceleration, with the largest gains at high frequency where the unaccelerated method is most expensive. Deflating the eigenvectors with eigenvalues closest to the driving frequency is the key, and choosing the number of deflation vectors to grow quadratically (in two dimensions) with $\omega$ keeps the asymptotic convergence rate roughly fixed as the frequency increases. On overset grids, where the discrete eigenvectors are accurate only to truncation error and are not orthogonal to machine precision, AUKED (AGMRES, ABICGSTAB) produced more accurate solutions than DEVD, while DEVD still gave results consistent with the resolution used; coarse-grid eigenvectors worked well for deflation provided they were adequately resolved. SVD compression of the eigenvectors reduced storage with no measurable effect on convergence, and the compressibility improved on finer and smoother (curvilinear) meshes. Together with the efficient EigenWave eigensolver, these results indicate that eigenvector deflation makes WaveHoltz an attractive solver for problems requiring many right-hand sides.



\bibliographystyle{elsart-num}
\bibliography{journal-ISI,ref,henshaw,henshawPapers,deflation}

\end{document}